\documentclass[12 pt]{amsart}
\usepackage[margin=1in]{geometry}
\usepackage[toc,page]{appendix}
\usepackage[notcite,notref]{} 
\usepackage{amsmath}
\usepackage{amsfonts}
\usepackage{amssymb}
\usepackage{amsthm}
\usepackage{setspace}
\usepackage{amsrefs}
\usepackage{inputenc}
\usepackage[english]{babel} 
\usepackage{comment}
\usepackage{enumitem}
\numberwithin{equation}{section} 

\usepackage{hyperref}
\hypersetup{hidelinks}

\title[]{Weighted Sobolev $L^{p}$ estimates for homotopy operators on strictly pseudoconvex domains with $C^{2}$ boundary}
\author[]{Ziming Shi}

\address{Department of Mathematics,
	University of Wisconsin-Madison, Madison, WI 53706}
\email{zimingshi@math.wisc.edu}


\newcommand{\dist}{\operatorname{dist}}

\newtheorem{thm}{Theorem}[section]
\newtheorem{cor}{Corollary}[thm]
\newtheorem{prop}[thm]{Proposition}
\newtheorem{lemma}[thm]{Lemma}

\theoremstyle{definition}
\newtheorem{defn}[thm]{Definition}

\theoremstyle{remark}
\newtheorem*{rem}{Remark}
\newtheorem*{clm}{Claim}
\newtheorem*{ack}{Acknowledgment}

\renewcommand{\th}[1]{\begin{thm}\label{#1}}
	\renewcommand{\eth}{\end{thm}}
\newcommand{\co}[1]{\begin{cor}\label{#1}}
	\newcommand{\eco}{\end{cor}}
\renewcommand{\le}[1]{\begin{lemma}\label{#1}}
	\newcommand{\ele}{\end{lemma}}
\newcommand{\pr}[1]{\begin{prop}\label{#1}}
	\newcommand{\epr}{\end{prop}}

\newcommand{\df}[1]{\begin{defn}\label{#1}}
	\newcommand{\edf}{\end{defn}}

\newcommand{\mk}{\begin{rem}}
\newcommand{\emk}{\end{rem}}
\newcommand{\cl}{\begin{clm}}
\newcommand{\ecl}{\end{clm}} 
\newcommand{\ac}{\begin{ack}}
\newcommand{\eac}{\end{ack}} 

\newcommand{\ga}{\begin{gather}}
\newcommand{\ega}{\end{gather}}
\newcommand{\gan}{\begin{gather*}}
\newcommand{\egan}{\end{gather*}}
\newcommand{\al}{\begin{gngn}}
\newcommand{\eal}{\end{align}}
\newcommand{\aln}{\begin{align*}}
\newcommand{\ealn}{\end{align*}}
\newcommand{\eq}[1]{\begin{equation}\label{#1}}
\newcommand{\eeq}{\end{equation}}

\newcommand{\ra}{\longrightarrow}

\newcommand{\sm}{\setminus}

\newcommand{\pa}{\partial{}}

\newcommand{\na}{\nabla}

\newcommand{\DD}[2]{\frac{\partial #1}{\partial #2}}

\newcommand{\R}{\mathbb{R}}
\newcommand{\C}{\mathbb{C}}

\newcommand{\ov}{\overline}
\newcommand{\ti}{\tilde}
\newcommand{\wti}{\widetilde}
\newcommand{\hht}{\widehat}

\renewcommand{\dbar}{\overline\partial}

\newcommand{\lip}{\operatorname{Lip}}

\newcommand{\mc}{\mathcal}

\newcommand{\all}{\alpha}

\newcommand{\del}{\delta}

\newcommand{\var}{\varphi}
\newcommand{\e}{\epsilon}
\newcommand{\ve}{\varepsilon}
\newcommand{\om}{\omega}
\newcommand{\Om}{\Omega}

\newcommand{\La}{\Lambda}
\newcommand{\la}{\lambda}

\newcommand{\gm}{\gamma}
\newcommand{\Gm}{\Gamma}

\newcommand{\yh}{\frac{1}{2}}

\newcommand{\re}[1]{(\ref{#1})}

\newcommand{\rl}[1]{Lemma~\ref{#1}}

\newcommand{\rp}[1]{Proposition~\ref{#1}}
\newcommand{\rt}[1]{Theorem~\ref{#1}}
\newcommand{\rd}[1]{Definition~\ref{#1}}

\newcommand{\supp}{\operatorname{supp}}

\newcommand{\we}{\wedge}

\newcounter{pp}
\newcommand{\bpp}{\begin{list}{$\hspace{-1em}\alph{pp})$}{\usecounter{pp}}}
\newcommand{\epp}{\end{list}}

\newcounter{ppp}
\newcommand{\bppp}{\begin{list}{$\hspace{-1em}(\roman{ppp})$}{\usecounter{ppp}}}
\newcommand{\eppp}{\end{list}}

\def\beq{\begin{equation}}
\def\eeq{\end{equation}}

\begin{document}
\begin{abstract}
	We derive estimates in a weighted Sobolev space $W^{k,p}_{\mu}(D)$ for a homotopy operator on a bounded strictly pseudoconvex domain $D$ of $C^2$ boundary in ${\C}^n$. As a result, we show that given any $2n < p < \infty$, $k > 1$, $q \geq 1$, and a $\dbar$-closed $(0,q)$ form $\var$ of class $W^{k,p}(D)$, there exist a solution $u$ to $\dbar u = \var$ such that $u \in W^{k,p}_{\yh-\ve}(D)$ for any $\ve > 0$. If $k=1$, then we can take $p$ to be any value between $1$ and $\infty$. In other words, the solution gains almost $\yh$-derivative in a suitable sense. 
\end{abstract} 	
	
\thanks{Supported in part by NSF grant DMS-1500162}  	
\maketitle
\tableofcontents

\section{Introduction}

In this paper we prove a regularity result concerning the solution of $\dbar$-equation on a strictly pseudoconvex domain $D$ with respect to a weighted Sobolev norm, assuming the boundary $bD$ is $C^{2}$. We define the \emph{weighted Sobolev space} $W^{k, p}_{\mu} (D)$ for a bounded domain $D \subset \R^{N}$ to be the subspace of $W^{k,p} (D)$ with  norm
\eq{wsobdef}
	 \| u \|_{W^{k,p} _{\mu} (D)}  = \sum_{|\all| \leq k+1} \left( \int_{D} | \pa^{\all} u(x) | ^{p} d(x) ^{ (1-\mu) p } \, dx  \right)^{\frac{1}{p}}. 
\eeq 
Here $k$ is a non-negative integer, $1 \leq p < \infty$, $0< \mu < 1$, and $d(x) = \dist(x, b D)$. These are Banach spaces with the norm defined as above. The reader can refer to \cite{HT91} for some properties of $W^{k,p}_{\mu} (D)$. 

 We mention some brief history regarding the ``$\yh$-estimate" for $\dbar$ solution $u$ of $\dbar u = \var$, for a $ \dbar$-closed $(0,q)$ form $\var$ on bounded strictly pseudoconvex domains. Regarding Sobolev space estimates,  Greiner and Stein \cite{G-S77} showed that for $q=1$, Kohn's canonical solution $\dbar^{\ast} N \var$ is in  $\mc{L}^{p}_{k+ \yh}(D)$, if $\var$ in $\mc{L}^{p}_{k} (D)$, for $1 < p < \infty$, and any non-negative integer $k$. Here $\mc{L}_{k}^{p}(D)$ is the Bessel potential space, as defined in \cite{SE70}*{p.~135}. Chang \cite{CD89} extended this result for all $q \geq 1$. Both Greiner-Stein and Chang assume that $b D$ is smooth. 
 
On the Hölder estimate side of $\dbar$ solutions, Henkin and Romanov \cite{H-R71} first achieved the $C^{\yh}$ estimate of $\dbar$ solutions for continuous $(0,1)$ form $\var$. Siu \cite{SY74} proved the $C^{k+ \yh}$ estimate for $q = 1$ and $k \geq 1$. Lieb-Range \cite{L-R80} constructed a $\dbar$ solution operator $H_{q}$, $q \geq 1$ and proved the $C^{k+ \yh}$ estimate when the boundary is $C^{k+2}$. In both results of Siu and Lieb-Range, $\var$ is assumed to be $\dbar$ closed. When $bD$ is smooth, Greiner and Stein (for $q = 1$) \cite{G-S77} showed that Kohn's canonical solution is in $\La_{r + \yh}$ if $\var \in \La_{r}$, for all $r>0$. Here $\La_{r}$ stands for the Zygmund space, as defined in \cite{G-S77}*{p.~141}. Chang \cite{CD89} extended this result for any $q \geq 1$ on the Siegal upper-half space. 

Recently Gong \cite{GX18} derived a new homotopy formula (see \re{Hq} and \re{H0} below), 
\begin{gather} \label{hfintro1}
  \var = \dbar H_{q} \var + H_{q+1} \dbar \var, \quad q \geq 1, \\ \label{hfintro2}
  \var = H_0 \var + H_{1} \dbar \var, \quad q=0. 
\end{gather}
for a bounded strictly pseudoconvex domain in $\C^{n}$ with the minimal smoothness condition of $C^{2}$ boundary. 
He showed that for $q \geq 1$, $H_{q} \var$ is in $\La_{r + \yh}$ if $\var \in \La_{r}$, $r>1$, and $H_{q} \var$ is in $C^{\frac{3}{2}}(\ov{D})$ if $\var \in C^{1}(\ov{D})$. Furthermore, the estimates do not require $\var$ to be $\dbar$-closed. 
There are two main features in the above homotopy formula in \cite{GX18}. The first is the regularized Leray map, introduced in \cite{GX18}. The second feature is the commutator $[\dbar, E]$, where $E$ is an extension operator bounded in $\La_{r}$-norm. This commutator was introduced by Peters \cite{PK91} and it has been used by Michel \cite{MJ91}, Range \cite{RM92}, Michel-Shaw \cite{M-S99}, Alexandre \cite{AW06} and others.

We shall prove our estimates for the homotopy operator $H_q$ and $H_0$. In section 2 we prove that homotopy formulas \re{hfintro1} and \re{hfintro2} hold in the distribution sense if $\var, \dbar \var \in W^{1,1} (D)$; see \rp{hfdist}. The goal is to prove the following: 
\th{mtintro}
   Let $D \subset \C^{n}$ be a bounded strictly pseudoconvex domain with $C^{2}$ boundary. Let $k$ be a positive integer, and $q$ be a non-negative integer. \\
   (i) Let $1 < p < \infty$, and $q > 0$. Then for any $\beta$, $0 < \beta < \yh$, 
   \[ 
     \| H_{q} \var \|_{W^{1,p}_{\beta}(D)} \leq C \| \var \|_{ W^{1,p}(D)}. 
   \] 
   (ii) Let $2n < p < \infty $, $k \geq 2$, and $q > 0$. Then  for any $\beta$, $0 < \beta < \yh$, 
   \[ 
      \| H_{q} \var \|_{W^{k,p}_{\beta}(D)} \leq C \| \var \|_{W^{k,p}(D)}.
   \]
   (iii) Let $1 < p < \infty$. Then for any $\beta$,  $0 < \beta < 1$, 
   \[ 
     \| H_{0} \var \|_{W^{0,p}_{\beta} (D)} \leq C \| \var \|_{W^{1,p}(D)}. 
   \] 
   (iv) Let $2n < p < \infty$, $k \geq 2$. Then for any $\beta$, $0 < \beta < 1$, 
   \[ 
     \| H_{0} \var \|_{W^{k-1,p}_{\beta} (D)} \leq C \| \var \|_{W^{k,p}(D)}.
   \]
Here we denote by $C$ some positive constants which depend on $D$, $n$, $p$, and $\beta$. 
\eth

We emphasize that $\var$ in the above estimates are not necessarily $\dbar$-closed. As a consequence, we have the following corollary: 
\co{}
   Let $D \subset \C^{n}$ be a bounded strictly pseudoconvex domain with $C^{2}$ boundary. Let $q$ be a positive integer. There exist a solution operator $H_{q}$ to the $\dbar$-equation $\dbar u = \var$ in $D$, for a given $\dbar$-closed $(0,q)$ form $\var$, such that the estimates in (i) and (ii) of \rt{mtintro} hold. In other words the solution $u$ gains $``\yh - \ve"$ derivative.
\eco

The paper is organized as follows. In section 2 we collect a few facts about the Stein extension operator, Sobolev space and the trace operator. We then derive the homotopy formula for Sobolev classes. We also recall from \cite{GX18} the regularized Leray map and its properties. 
In section 3 we prove the estimates for $H_q$, $q \geq 1$ (part (i) and (ii) of \rt{mtintro}). The main techinical part involves a subtle use of integration by parts to move derivatives from the kernel to $\var$. In section 4 we prove the estimates for $H_0$ (part (iii) and (iv) of \rt{mtintro}). 

\ac{} 
   I am grateful to my advisor Professor Xianghong Gong for his guidance. 
\eac

\section{Homotopy formula for Sobolev Space}
In this section we derive the homotopy formula introduced in \cite{GX18} for the Sobolev classes. We shall need some standard facts about Sobolev spaces. For reader's convenience we state them here. We use $W^{k,p}(D)$ to denote the usual Sobolev space with norm 
\[
    \| u \|_{W^{k,p}(D)}  = \sum_{|\all| \leq k} \left( \int_{D} | \pa^{\all} u(x) | ^{p} \, dx  \right)^{\frac{1}{p}}. 
\]
We remind the reader that the $\dbar$ solution space $W^{k,p}_{\mu} (D)$ defined in section 1 has actually $k+1$ interior derivatives. Thus $W^{k,p}_{\mu} (D) \subset W^{k+1,p} (D')$, for any relatively compact subdomain $D'$ of $D$. 
\pr{Sobemb}
Let $D \subset \R^{N}$ be a bounded domain with $C^{1}$ boundary. Assume $N < p \leq \infty$ and $u \in W^{k,p} (D)$. Then up to a set of measure $0$, $u \in C^{k-1, \all}(\ov{D})$, for $\all = 1 - \frac{N}{p} > 0$, and $u$ satisfies the estimate 
\[ 
   \| u \|_{C^{k-1,\all} (\ov{D})} \leq C \| u \|_{W^{k,p}(D)}, 
\] 
where $C$ depends on $k,p,N$ and $D$. 
\epr

The proof can be found in \cite{LG09}*{p.~335}.

We need an extension operator due to E. Stein. 
\pr{extthm}
   Let $D$ be a bounded domain whose boundary satisfies the minimal smoothness condition as defined in \cite{SE70}*{p.~189}, (in particular, a bounded domain is minimally smooth if its boundary is locally given by graphs of Liptschitz functions.) Then 
   (i) There is a continuous linear operator $E: W^{k,p}(D) \to W^{k,p}(\R^{N}) $ so that $Ef = f$ on $D$, for all $p$, $1 \leq p \leq \infty$, and all non-negative integer $k$.  \\
   (ii) There is a continuous linear operator $E: C^{0}(\ov{D}) \to C^{0} (\R^{N})$ so that $Ef = f$ on $D$ and 
   \[
     | Ef |_{C^{r} (\R^{N})} \leq C(r,D) | f |_{C^{r} (\ov{D})}, \quad \forall \, r \in [ 0, \infty).   
   \]
\epr 
The proof of (i) can be found in \cite{SE 70}*{p.~181}, and the proof of (ii) can be found in \cite{GX18}.  

In what follows we denote $\R^{N}_{+} = \{x = (x', x_N) \in \R^{N-1} \times \R, x_N > 0 \}$. 
\df{unfLip}
  The boundary $b \om$ of an open set $\om\subset \R^{N}$ is \emph{uniformly Lipschitz} if there exist $\ve, L > 0$, $M \in \mathbb{N}$, and a locally finite countable open cover $U_{l}$ of $b \om$ such that \\
  (i) If $x \in b \om$, then $B(x, \ve) \subset U_{l} $ for some $l \in \mathbb{N}$. \\
  (ii) No point of $\R^{N}$ is contained in more than $M$ of the $U_l$'s. \\
  (iii) For each $k$ there exist local coordinates $y = (y_1, \dots, y_{\all_l},  \dots, y_N)$ and a Lipschitz function $f_l: \R^{N-1} \to \R$ with $\lip f_l \leq L$, such that
  \eq{Al}
    U_l \cap \om = U_l \cap A_l, \quad A_l = \{ y \in \R^{N}: y_{\all_l} > f_l(y'_{\all_l}) \}. 
  \eeq
  where $y'_{\all_l} = (y_1, \dots, \hht{y_{\all_l}}, \dots, y_N)$, and $\hht{\cdot}$ means $\cdot$ is omitted. 
\edf
We now define the trace operator for $W^{1,1} (\om)$. First we define it on $W^{1,1} (\R^{N}_{+})$. 
\pr{}
Let $N \geq 2$ and let $W^{1,1}_0(\R^{N}_{+})$ be the family of all functions $u \in W^{1,1} (\R^{N}_{+})$ with bounded support. Then there exist a linear operator 
\[ 
  Tr: W^{1,1}_0(\R^{N}_{+}) \ra L^{1} (\R^{N-1}) 
\]
such that \\
(i) $Tr (u)(x') = u (x', 0)$ for all $x' \in \R^{N-1}$, and for all $u \in W^{1,1}_0(\R^{N}_{+}) \cap C (\ov{\R^{N}_{+}})$.  \\
(ii) For all $u \in W^{1,1}_0(\R^{N}_{+})$, 
\eq{tresths}
\int_{\R^{N-1}} \left| Tr(u) (x') \right| \, dx' \leq \int_{\R^{N}_{+}} \left| \DD{u}{x_N} (x) \right| \, dx. 
\eeq 
(iii) For all $v \in C^{1}_{c}(\R^{N})$, $u \in W^{1,1}_0(\R^{N}_{+})$, and $i = 1, \dots, N$, 
\eq{stokeshs}
\int_{\R^{N}_{+}} \DD{(uv)}{x_i} \, dx = \int_{\R^{N-1}} u Tr(v) \nu_i \, dx'
\eeq
where $\nu = - e_N = (0, \dots, 0, -1)$ is the outer unit normal on $\R^{N-1} = \{ x_N = 0 \}$, $dx = dx_1 \dots dx_N$ and $dx' = dx_1 \dots d x_{N-1}$.  
\epr
For proof see \cite{LG09}*{p.~452}.
\pr{trthm} 
 Let $\om \subset \R^{N}$, $N \geq 2$, be an open set whose boundary $b \om$ is uniformly Lipschitz, with the corresponding $\ve, L, M$ given as in \rd{unfLip}. There exist a continuous linear operator 
\[
  Tr: W^{1,1} (\om) \ra L^{1} (b \om, ds_{b \om}) 
\]
such that \\
(i) $Tr(u) = u$ on $b \om$ for all $u \in W^{1,1} (\om) \cap C(\ov{\om})$. \\
(ii) Denote by $ds_{b \om}$ the surface element of $b \om$. We have
\eq{trest} 
  \int_{b \om} \left| Tr(u) \right| \, ds_{b \om} \leq \frac{CM}{\ve} \sqrt{1+ L^{2}} \int_{\om} |u| \, dx + \sqrt{1+ L^{2}} \int_{\om} | \na u | \, dx.  
\eeq
\epr
The reader can refer to \cite{LG09}*{p.~460-462} for the proof of \rp{trthm}. For later use we recall the construction of the above trace operator. Let $\{U_{l} \}$ be an open cover of $b \om$ as given in \rd{unfLip}, and let $\chi_l$ be smooth partition of unity such that $\supp \chi_l \subset \subset U_I$. Then $u = \sum_{l} \chi_l u := \sum_{l} u_{l} $ in a neighborhood of $b \om$, and $u_l$ has compact support in $U_l$. Since $\om \cap U_l = A_l \cap U_l$ (\re{Al}), we can extend $u_{l}$ to be $0$ in $A_l \sm U_l$ to obtain $u_l \in W^{1,1} (A_l)$. Define
\eq{trdeffcn}
Tr(u) = \sum_l Tr(u_l), \quad Tr(u_l) = Tr(u_l \circ \psi_l) \circ \psi_l^{-1}, \quad u_l := \chi_l u. 
\eeq
where $\psi_l: \R^{N} \to \R^{N}$ is given by $\psi_l(y) = (y_1, y_{\all_{l-1}}, y_{\all_l} + f_l(y_{\all_l}'), y_{\all+1}, \dots, y_{N}  )$. 
Furthermore we can choose the partition of unity $\chi_{l}$ so that $Tr(u_l)$ is compactly supported in $b \om \cap U_{l}$.   

Let $\phi = \sum_{I} \phi^{I} dx_I$ be a differential form of degree $q$, for $q \geq 1$. We say that $\phi \in W^{1,1}(\om)$ if each component function $\phi^{I}$ belongs to the class $W^{1,1} (\om)$. We define the trace of $\phi$ on $b \om$ to be
\eq{trdefdf}
  Tr(\phi) = \sum_{| I | =q } Tr(\phi^{I}) dx_{I}. 
\eeq

\pr{stokesdist}
Let $\om \subset \R^{N}$ be a bounded domain with uniformly Lipschitz boundary. Suppose that $\phi$ is a differential form and $\phi \in W^{1,1} (\om)$. We have
\eq{stokeseqn}
\int_{b \om} Tr(\phi) \we \all = \int_{\om} d (\phi \we \all) 
\eeq
for any $\all$ which is a $C^{1}_{0}(\R^{N})$ form. 
\epr
Formula \re{stokeseqn} can be proved by pulling back the forms to the upper half plane $\R^{N}_{+}$ by Lipschitz maps, smoothing out the Lipschitz maps and using \re{stokeshs}. We leave the details to the reader.   
   
\le{trbdycvg}
Let $\om \subset \R^{N}$ be a bounded domain with $C^{1}$ boundary, and $\om' \subset \subset \om $. Suppose $k(z, \zeta)$ is uniformly bounded for $z \in \om'$ and $\zeta$ in some neighborhood of $b \om$, and is uniformly continuous in $\zeta$. Suppose $u \in W^{1,1} (\om)$. Let $\om_{j}$ be a sequence of smooth domains approximating $\om$ from inside, i.e. $\om_{j} \subset \subset \om_{j+1} \subset \subset \cdots \om$, and such that locally the defining functions of $b \om_j$ converge uniformly to that of $b \om$ in $C^{1}$-norm, and  Then
\[ 
  \int_{b \om_{j}}  k(z, \zeta)  Tr(u)(\zeta) \,  ds_{b \om_{j}}(\zeta)  \overset{j \to \infty}{\ra}\int_{b \om}  k(z, \zeta)  Tr(u)(\zeta) \,  ds_{b \om} (\zeta)
\]
uniformly on $z \in \om'$. Here $ds_{b \om} (\zeta)$ and $ds_{b \om_j}(\zeta)$ denote the surface elements of $b \om$ and $b \om_j$ respectively.  
\ele
\begin{proof}
	Let $\{ U_l \} $ be a (finite) open cover of $b \om$ and $b \om_j$ as given in \rd{unfLip}, for $j$ sufficiently large. By the way we define trace \re{trdeffcn}, it suffices to prove that for each $l$, 
\[
   \int_{b \om_j \cap U_{l} }  k(z, \zeta)  Tr(u)(\zeta) \,  ds_{b \om_{j}}(\zeta)  \overset{j \to \infty}{\ra}\int_{b \om \cap U_l}  k(z, \zeta)  Tr(u)(\zeta) \,  ds_{b \om} (\zeta)
\]
where $u$ has compact support in $U_l$. There exist local coordinates $x = (x', x_N) \in \R^{N-1} \times \R$, and $C^{1}$ functions $f$, $f_j$, $\R^{N-1} \to \R$, such that $ \om \cap U_l = A \cap U_l$, and $\om_j \cap U_l = A_j \cap U_l$, where
\[ 
  A = \{ (x', x_{N}) \subset \R^{N}: x_{N}  > f(x') \}, \quad A_j = \{ (x', x_{N}) \subset \R^{N}: x_{N}  > f_{j}(x') \}. 
\]  
Since $\om_j \subset \subset \om$, we can assume $A_j \subset \subset A$. 
Since $u$ has compact support in $U_{l}$, we can extend $u$ to be $0$ in $A \sm U_{l}$ (Thus also $0$ in $A_j \sm U_l$.) to obtain $u \in W^{1,1} (A)$ and $u \in W^{1,1} (A_j)$. By assumption, $f_{j}$ converges uniformly to $f$ in $ C^{1} (\R^{N-1})$. The surface area element on $b \om \cap U_l$ is given by
\[
   ds(b \om) = \sqrt{1+ \left| \na f (x') \right|^{2} } dx', \quad dx' =  dx_{1} \cdots dx_{N-1}, 
\]
and similarly $ds(b \om_{j}) = \sqrt{1+ \left| \na f_{j} (x') \right|^{2} } dx'$. Define $C^{1}$ diffeomorphisms $\psi, \psi_j : B_0 \to U_0$ by 
\gan
  \psi (x) = (x', x_{N} + f(x')),
\\
  \psi_{j} (x ) = (x', x_{N} + f_{j}(x')). 
\end{gather*}
Let $\wti{u} = u \circ \psi$ , $\wti{u}_j = u \circ \psi_j$. 
	Note that $\psi: \R^{N}_{+} \to A$, and $\psi_j : \R^{N}_{+} \to A_j$, and $\wti{u}$ and $\wti{u_{j}}$ are functions in $W^{1,1}_0(\R^{N}_{+})$. By \re{trdeffcn}, $Tr(u)|_{b A} (x', f(x')) = Tr(\wti{u}) (x')$ and 
	$Tr(u)|_{bA_j} (x', f_{j}(x')) = Tr(\wti{u}_{j}) (x' )$, for $x' \in \R^{N-1}$. As remarked before, $Tr(u)|_{bA}$ (resp. $Tr(u)|_{bA_j}$) is compactly supported in $b A \cap U_l$ (resp. $b A_j \cap U_l$). Since
	\[ 
	  Tr(\wti{u}) = Tr(u \circ \psi) = Tr \left. \right|_{bA} (u) \circ \psi, \quad Tr(\wti{u}_j) = Tr(u \circ \psi_j) = Tr(u)|_{bA_j} \circ \psi_j, 
	\]
    $Tr(\wti{u})$ and $Tr(\wti{u}_j)$ are compactly supported in $\R^{N-1}$. 
	Thus
\begin{align} \label{FGH}
	&\left| \int_{b \om_{j} \cap U_l}  k(z, \zeta)  Tr(u)(\zeta)  ds(b \om_{j}) -  \int_{b \om \cap U_l}  k(z, \zeta)  Tr(u)(\zeta)  ds(b \om) \right| 
	\\ \nonumber & = \left| \int_{\R^{N-1}} k(z, (x', f_j(x'))) Tr(\wti{u_{j}}) (x') g_{j}(x') - k (z, (x', f(x'))) Tr(\wti{u}) (x') g(x')  dx' \right| 
	\\ \nonumber & \quad \leq  F_{j} (z) + G_{j}(z) + H_j(z), 
\end{align}
where 
\[
  F_{j} (z) = \int_{\R^{N-1}}  \left| k(z, (x', f_j(x'))) - k(z, (x', f(x'))) \right|  \left| Tr(\wti{u_{j}}) (x') g_{j}(x') \right| \, dx', 
\]
\[ 
  G_j (z) =  \int_{\R^{N-1}}  \left| k(z, (x', f(x'))) \right|  \left| Tr(\wti{u_{j }} - \wti{u})  (x') \right| |g_j (x')|  \, dx', 
\]
\[
  H_j (z) =  \int_{\R^{N-1}}  \left| k(z, (x', f(x'))) \right|  \left| Tr(\wti{u})  (x') \right| | (g_j - g )(x')|  \, dx', 
\]
and
\[ 
  g(x') = \sqrt{1+ \left| \na f (x') \right|^{2} }, \quad g_j (x')= \sqrt{1+ \left| \na f_j (x') \right|^{2} }. 
\]
By assumption, $| f_j - f |$ converges to $0$ uniformly on $\R^{N-1}$ and $k(z, \zeta)$ is uniformly continuous in $\zeta$ in a neighborhood of $b \om$, so we have
\[ 
\left| k(z, (x', f_j(x'))) - k(z, (x', f(x'))) \right| \overset{j \to \infty}{\ra} 0, \quad \text{uniformly in $x' \in \R^{N-1}$}. 
\]
Hence to show $F_j$ converges to $0$ uniformly in $z \in \om'$, it suffices to show 
\eq{trbdd}
  \int_{B_1^{N-1}}  \left| Tr(\wti{u_{j}}) (x') \right| \left| g_{j}(x') \right| \, dx'  \leq C
\eeq
for some $C$ independent of $j$. By \re{tresths}, we have
\begin{align*}
\int_{\R^{N-1}}  \left| Tr(\wti{u_{j}}) (x') g_{j}(x') \right| \, dx' 
&\leq C \int_{\R^{N}_{+}} \left| \DD{\wti{u_j}}{x_N} (x) \right| \, dx \\
&\leq \int_{\R^{N}_{+}} \left| \DD{u}{x_N} (x', x_N + f_j(x')) - \DD{u}{x_N} (x', x_N + f(x')) \right| \, dx \\
&\quad +  \int_{\R^{N}_{+}} \left| \DD{u}{x_N} (x', x_N + f(x')) \right| \, dx. 
\end{align*}
Below we show the first integral in the last inequality converges to $0$ as $j \to \infty$. This proves \re{trbdd} and thus $F_j$ converges to $0$ uniformly in $z \in \om'$. By assumption, $ |g_{j} - g|$ converges to $0$ uniformly on $\R^{N-1}$, $|k(z, (x', f(x'))) | \leq C$ for $x' \in \R^{N-1}$, and $Tr(\wti{u}) \in L^{1} (\R^{N-1})$,  it follows that $H_j$ converges to $0$ uniformly on $z \in \om'$. For $G_j$, by \re{tresths} we have
\begin{align} \label{Fjest}
\int_{\R^{N-1}} \left| Tr(\wti{u} - \wti{u_j}) (x') \right| \, dx' 
& \leq \int_{\R^{N}_{+}} \left| \DD{\wti{u}}{x_N} (x) 
- \DD{\wti{u_j}}{x_N} (x) \right| \, dx 
\\ \nonumber & =  \int_{\R^{N}_{+}} \left| \DD{u}{y_N} (\psi(x)) - \DD{u}{y_N} (\psi_j (x))\right| \, dx. \\ \nonumber 
&= \int_{A} \left| \DD{u}{y_N} (y) - \DD{u}{y_N} (y', y_N + (f_j - f) (y')) \right| \, dy. 
\end{align}  
Since $f_j$ converges to $f$ uniformly on $\R^{N-1}$, we can show the last integral converges to $0$ by a standard smoothing argument. Since
\[ 
\left| k(z, (x', f(x'))) \right| \leq C, \quad |g_j (x')| \leq C, 
\]
we have proved that $G_{j}$ converges to $0$ uniformly on $ z \in \om'$. The conclusion of the lemma then follows from estimate \re{FGH}.  
\end{proof}

We now extend the homotopy formula in \cite{GX18} to $\var$ satisfying $\var, \dbar \var \in W^{1,1}(D)$. 
\pr{hfdist}
   Let $D \subset \C^{n}$ be a bounded domain with $C^{1}$ boundary and let $U$ be a bounded neighborhood of $\ov{D}$. Let $g^{0} = \ov{\zeta -z}$. Let $g^{1} = W(z, \zeta)$, where $W \in C^{1}(D \times (U \sm D))$ is a Leray mapping, that is, $W$ is holomorphic in $z \in D$ and satisfies 
   \[ 
     \Phi(z, \zeta) = W(z, \zeta) \cdot (\zeta -z) \neq 0, \quad z \in D, \quad \zeta \in U \sm D. 
   \] 
   Let $\var$ be a $(0,q)$-form. Suppose that $\var$ and $\dbar \var$ are in $W^{1,1}(D)$. (That is all the coefficient functions of $\var$ and $ \dbar \var$ are in $W^{1,1} (D)$). Then we have the following: \\
(i) The Bochner-Martinelli formula 
\eq{propBMKf} 
   \var = \dbar_{z} \int_{D} \Om^{0}_{0, q-1} (z, \zeta) \we \var + \int_{D} \Om^{0}_{0,q} (z, \zeta) \we \dbar \var  
  + \int_{b D} \Om^{0}_{0, q} (z, \zeta) \we Tr(\var)
\eeq 
holds in the distribution sense in $D$. 
\\
(ii) The following homotopy formula holds in $D$ in the distribution sense. 
\eq{homf1}
    \var = \dbar H_{q} \var + H_{q+1} \dbar \var, \quad 1 \leq q \leq n 
\eeq
\eq{homf2}
   \var = H_{0} \var + H_{1} \dbar \var, \quad q = 0 
\eeq
where
\eq{Hq}
   H_{q} \var := \int_{U} \Om^{0}_{0, q-1} \we E \var + \int_{U \sm D} \Om^{0,W}_{0, q-1} \we [\dbar, E] \var, \quad 1 \leq q \leq n
\eeq
\eq{H0}
   H_{0} \var := \int_{U \sm D} \Om_{0,0}^{1}  \we [\dbar, E] \var , \quad \quad [\dbar, E] \var = \dbar E \var - E \dbar \var. 
\eeq
Here $\Om^{\bullet}_{0,q}$ stands for the $(0,q)$ component of $\Om^{\bullet}$ of type $(0, q)$ in $z$, and 
\eq{Om0}
  \Om^{0} (z, \zeta) = \frac{1}{(2 \pi i)^{n}} \frac{\left<\ov{\zeta} - \ov{z} \, , \, d \zeta \right>}{|\zeta -z |^{2}} \we \left( \dbar_{\zeta,z} \frac{ \left< \ov{\zeta} - \ov{z} \, , \, d \zeta \right>}{|\zeta -z|^{2}} \right)^{n-1}, \quad \dbar_{\zeta,z} = \dbar_{\zeta} + \dbar_{z}; 
\eeq
\eq{Om1}
   \Om^{W} (z, \zeta) = \frac{1}{(2 \pi i)^{n}} \frac{\left< W, d \zeta \right>}{\Phi (z, \zeta)} \we \left[ \dbar_{\zeta, z} \frac{ \left< W, d \zeta \right>}{\Phi(z,\zeta)} \right]^{n-1}, \quad \Phi (z, \zeta) = W(z, \zeta) \cdot (\zeta -z); 
\eeq
\begin{align} \label{Om01}
   &\Om^{0,W} (z, \zeta) 
   = \frac{1}{(2 \pi i)^{n}} \frac{\left<\ov{\zeta} - \ov{z} \, , \, d \zeta \right>}{|\zeta -z |^{2}} \we \frac{\left< W, d \zeta  \right>}{\left< W \, , \, \zeta -z \right>} 
       \\ \nonumber & \qquad \we \sum_{i+j=n-2} \left[ \frac{\left< d\ov{\zeta} - d\ov{z} \, , \, d \zeta \right>}{|\zeta -z |^{2}} \right]^{i}   \we \left[ \dbar_{\zeta,z} \frac{ \left< W  ,  d \zeta  \right>}{\left< W, \zeta -z \right>} \right]^{j}.
\end{align}
We set $\Om^{W}_{0,-1} = 0$ and $\Om^{0,W}_{0,-1} =0$. 
\epr
\begin{proof}  
(i) By some abuse of notation, we shall denote the coefficeint functions of $\var$ by $\var$, and our smoothing is done componentwise. Let $\{ \psi_{\ve} \}_{\ve >0}$ be the standard mollifier which satisfies $\psi_{\ve} \in C^{\infty}_{0} (B_{\ve}(0))$, $\psi_{\ve} \geq 0$, and $\int_{\C^{n}} \psi_{\ve} = 1$. Let $\var_{\ve} =  \var \ast \psi_{\ve} $ be defined by    
\gan
\var_{\ve} (z) = (z) =  \int_{\C^{n}} \var (z - \zeta) \psi_{\ve} (\zeta) \, d V(\zeta)   
= \int_{B_{\ve} (0)} \var (z - \zeta) \psi_{\ve} (\zeta) \, d V(\zeta). 
\end{gather*}
Then we can show that for any $D' \subset \subset D$, and $\ve < \ve_0$ sufficiently small, 
\eq{varcvgDj}
  \var_{\ve} \overset{\ve \to 0}{\ra} \var \quad  \text{in $W^{1,1} (D')$}, \quad  \dbar \var_{\ve} \overset{\ve \to 0}{\ra} \dbar \var \quad \text{in $W^{1,1} (D')$}. 
\eeq

When $bD \in C^{1}$ and $\var \in C^{1}(\ov{D})$, the proof of formula \re{propBMKf} can be found in \cite{C-S01}*{p.~265}. Let $D_{j}$ be a sequence of domains with $C^{\infty}$ boundary approximating $D$ from inside, $D_{j} \subset \subset D_{j+1} \subset \subset \cdots D$, and locally the defining functions of $D_j$ converge uniformly in $C^{1}$-norm. Fix $j$ and $\ve_{0}>0$ such that $dist(D_{j}, D) > \ve_{0}$. The formula \re{propBMKf} then holds for $\var_{\ve}$ on $D_{j}$, for any $\ve < \ve_{0}$: 
\[ 
\var_{\ve}(z) = \dbar_{z} \int_{D_{j}} \Om^{0}_{0, q-1} (z, \zeta) \we \var_{\ve} + \int_{D_{j}} \Om^{0}_{0,q} (z, \zeta) \we \dbar \var_{\ve} 
+ \int_{b D_{j}} \Om^{0}_{0, q} (z, \zeta) \we \var_{\ve}. 
\]  
By Sobolev embedding \cite{LG09}*{p.~312}, $W^{1,1} (D) \subset L^{\frac{2n}{2n-1}}(D)$.  Applying this and the Calderón-Zygmund estimate for the Newtonian potential \cite{G-T01}*{p.~230}, we have for any $D' \subset \subset D_j$, 
\begin{align} \label{Bqvarcvg}
\left| \int_{D_{j}} \Om^{0}_{0, q-1} (z, \zeta) \we (\var_{\ve} - \var) \right|_{W^{1,\frac{2n}{2n-1}}(D')} 
&\leq C(n) \| \var_{\ve} - \var \|_{L^{\frac{2n}{2n-1}}(D_{j})} \\
& \nonumber \leq C(n) \| \var_{\ve} - \var \|_{W^{1,1}(D_{j})}. 
\end{align} 
and similarly, 
\begin{align}\label{Bqdvarcvg}
\left| \int_{D_{j}} \Om^{0}_{0, q} (z, \zeta) \we (\dbar \var_{\ve} - \dbar \var) \right|_{W^{1,\frac{2n}{2n-1}}(D')} 
\leq  C(n) \| \dbar \var_{\ve} - \dbar \var \|_{W^{1,1}(D_{j})}. 
\end{align}
Note that the above constants $C$ depend only on the dimension $n$ and is independent of $j$. Now, $ | \Om^{0}_{0,q} (z, \zeta) | \leq C$ for $z \in D'$ and $\zeta \in bD_j$, $D' \subset \subset D_j$. By estimate \re{trest}, there exist a constant $C$ independent of $j$ such that
\begin{align} \label{Bqbdycvg} 
\left| \int_{b D_{j}} \Om^{0}_{0, q} (z, \zeta) \we (Tr(\var_{\ve}) - Tr(\var)) \right|_{C^{0}(\ov{D'})} 
&\leq C \| Tr(\var_{\ve}) - Tr(\var) \|_{L^{1}(b D_{j})} \\
& \nonumber \leq  C\| \var_{\ve} - \var \|_{W^{1,1}(D_{j})}. 
\end{align}
As $\ve \to 0$, all these expressions in \re{Bqvarcvg}, \re{Bqdvarcvg} and \re{Bqbdycvg} converge to $0$. Thus
\[ 
\var(z) = \dbar_{z} \int_{D_{j}} \Om^{0}_{0, q-1} (z, \zeta) \we \var + \int_{D_{j}} \Om^{0}_{0,q} (z, \zeta) \we \dbar \var 
+ \int_{b D_{j}} \Om^{0}_{0, q} (z, \zeta) \we \var. 
\] 
holds in the distribution sense. (In fact, we only need to show convergence in $L^{1}(D')$.) Finally we let $j \to \infty$. For some constant $C$ independent of $j$, we have
\[
\left| \int_{D \sm D_{j}} \Om^{0}_{0,q-1}(z,\zeta) \we \var \right| \leq C \| \var \|_{L^{1}(D \sm D_{j})}
\overset{j \to \infty}{\ra} 0, 
\]
\[ 
\left| \int_{D \sm D_{j}} \Om^{0}_{0,q} (z, \zeta) \we \dbar \var \right| \leq C \| \dbar \var \|_{L^{1}(D \sm D_{j})} \overset{j \to \infty}{\ra} 0, 
\] 
where the convergence is uniform on $D'$. For $z \in D'$, and  $\zeta$ in a small neighborhood of $bD$, $\Om^{0}_{0,q}(z, \zeta)$ is smooth in both variable. Hence we can apply \rl{trbdycvg} to get 
\eq{BMKfbie}
   \int_{b D_{j}} \Om^{0}_{0,q}(z, \zeta) \we Tr(\var) \overset{j \to \infty}{\ra} \int_{b D} \Om^{0}_{0,q}(z, \zeta) \we Tr(\var). 	
\eeq
Consequently the Bochner-Martinelli formula \re{propBMKf} holds in the distribution sense for any $\var$ satisfying $\var, \dbar \var \in W^{1,1}(D)$. 
\\ \\
(ii) We prove formula \re{homf1}. The proof for \re{homf2} is similar and we shall omit the proof. First let us derive \re{homf1} under the assumption $bD \in C^{2}$, $W\in C^2(D\times(U\setminus D))$ and $\var, \dbar \var \in C^{1}( \ov{D})$. This part of the proof is the same as presented in \cite{GX18}, and we put it here since later on we shall prove the same thing under weaker assumptions. The Bochner-Martinelli holds in this case: 
\eq{BMKf}
   \var = \dbar_{z} \int_{D} \Om^{0}_{0, q-1} (z, \zeta) \we \var + \int_{D} \Om^{0}_{0,q} (z, \zeta) \we \dbar \var 
   + \int_{b D} \Om^{0}_{0, q} (z, \zeta) \we \var. 
\eeq
For $q\geq1$, 
\eq{homiden}
   \Om_{0,q}^0 - \Om_{0,q}^{W} =\dbar_\zeta\Omega^{0,W}_{0,q}+\dbar_z\Omega^{0,W}_{0,q-1}, \quad (z, \zeta)\in  D \times (U \sm D).  
\eeq 
For the proof of this identity the reader can refer to \cite{C-S01}*{p.~264}. Applying this to the boundary integral in \re{BMKf} we get
\begin{align}\label{hfibp0}
\varphi(z)&=\dbar B_{q} \varphi(z)+B_{q+1} \dbar\varphi(z) + \int_{bD} \Om^{W}_{0,q} (z, \zeta) \we \var \\
&\quad\nonumber
+ \dbar_{z} \int_{b D}\Omega_{0,q-1}^{0,W}(\zeta,z)\wedge\varphi(\zeta)
- \int_{b D}\Omega_{0,q}^{0,W}(\zeta,z)\wedge\dbar_\zeta\varphi(\zeta), \quad z\in D, 
\end{align}
where we denote
\[ 
  B_{q} \var = \int_{D} \Om^{0}_{0, q-1}(z, \zeta) \we \var. 
\] 
We denote by $B_{q; \, \Om} (\var)$ for the above integral when the domain of integration is $\Om$, and $B_{q} \var = B_{q;\,D} (\var)$. Since $W$ is a Leray map and it is holomorphic in $z$, in view of expression \re{Om1},  
\eq{Om1=0}  
  \Om^{W}_{0,q} (z, \zeta) = 0, \quad \quad \text{for $q \geq 1$.}
\eeq
For the last two integrals in \re{hfibp0} we first extend $\var, \dbar \var$ to $E\var, E \dbar \var \in W^{1,1} (U)$ by means of \rp{extthm}. Applying Stokes theorem to the domain $U \sm D$ we get 
\begin{align} \label{hfibp1}
   \int_{b D}\Omega_{0,q-1}^{0,W} \wedge\varphi 
   &= \int_{U \sm D} \dbar_{\zeta} \Om^{0,W}_{0,q-1} \we E \var + \int_{U  \sm D}  \Om^{0,W}_{0, q-1} \we \dbar_{\zeta} E \var
   \\ \nonumber &= \int_{U \sm D}  \Om^{0} _{0,q-1} \we E \var-  \int_{U \sm D}  \Om^{W} _{0,q-1} \we E \var 
   \\ \nonumber &\qquad - \dbar_{z} \int_{U \sm D} \Om^{0,W}_{0, q-2} \we E \var + \int_{U  \sm D}  \Om^{0,W}_{0, q-1} \we \dbar_{\zeta} E \var, 
\end{align} 
and
\begin{align} \label{hfibp2} 
  &\int_{b D}\Omega_{0,q}^{0,W} \wedge E \dbar_\zeta \var
  = \int_{U \sm D} \dbar_\zeta \Om^{0,W}_{0,q} \we E\dbar_\zeta \var + \int_{U \sm D} \Om^{0,W}_{0,q} \we \dbar_{\zeta} E \dbar_{\zeta} \var  
  \\ \nonumber &\qquad = \int_{U \sm D} \Om^{0}_{0,q}  \we E\dbar_\zeta \var - \dbar_{z} \int_{U \sm D} \Om^{0,W}_{0, q-1} \we E \dbar_\zeta \var +\int_{U \sm D} \Om^{0,W}_{0,q} \we \dbar_{\zeta} E \dbar_{\zeta} \var. 
\end{align} 
By \re{Om1=0}, $\Om^{W}_{0,q-1} = 0 $ if $q \geq 2$. If $q = 1$, $\Om^{W}_{0, q-1} = \Om^{W}_{0,0}$ is holomorphic in $z$, and $\Om^{0,W}_{0, q-2} = \Om^{0,W}_{0, -1} =0$. Using these facts and subsituting \re{hfibp1} and \re{hfibp2} into \re{hfibp0}, we obtain \re{homf1}.  

Suppose now that $b D \in C^{1}$, $W\in C^1(D\times(U \setminus D))$ and $\var, \dbar \var \in W^{1,1}(D)$. We shall derive the homotopy formula \re{homf1} in the distribution sense. We need to justify \re{hfibp0}, \re{hfibp1} and \re{hfibp2}. 

As before, we take a sequence of domains $D_{j}$ with smooth boundary approximating $D$ from inside, such that locally the defining functions of $D_j$ converge in the $C^1$~norm. Consider sufficiently large $j$ such that $D' \subset \subset D_{j} \subset \subset D$. Let $\var_{\ve}$ be a sequence of smooth forms so that $\var_{\ve} \to \var$, and $\dbar \var_{\ve}  \to \dbar \var$ in $W^{1,1}(D_{j})$. Since $D \times (U \sm D)$ has $C^{1}$ boundary, by \rp{extthm} (ii) we can extend $W$ to get $EW \in C^{1} (\C^{n} \times \C^{n})$, such that $EW (z, \zeta) = W(z, \zeta)$ for $z \in D$ and $\zeta \in \ov{U \sm D} $. Note that $EW(\cdot, \zeta)$ may not be holomorphic for $\zeta \in D$. 

For $z \in D$ and $\zeta \in U$, define 
\eq{EWapprox}
(EW)_{\ve'}(z, \zeta)=\int_{\C^{n} \times \C^{n}} \psi_{\ve'} (z' - z, \zeta'-\zeta)EW(z', \zeta')\, dV (z') dV(\zeta'). 
\eeq
where $\psi_{\ve'}$ is the standard mollifier. Then $(EW)_{\ve'}$ is $C^\infty$ in $\C^{n} \times \C^{n}$. 
Also $\left< (EW)_{\ve'}, \zeta -z \right> \neq 0$ for $z \in D'$ and $\zeta \in b D_j$, if $\ve'$ is sufficiently small and $j$ is sufficiently large. Indeed, by assumption $\left<W, \zeta -z \right> \neq 0$ on $D \times (U \sm D)$. Since $\ov{D'} \times bD$ is a compact subset of $D \times (U \sm D)$, $\left| \left< EW, \zeta -z \right> \right| = \left| \left< W, \zeta -z \right> \right| \geq \del$ on $D' \times b D$, and $\left< (EW)_{\ve'}, \zeta -z \right> \geq \del' $ if $\ve'$ is small and $j$ is large. 
\begin{align} \label{Omapprox}
   \Om^{0,\, (EW)_{\ve'}} (z, \zeta) &= \frac{1}{(2 \pi i)^{n}} \frac{\left<\ov{\zeta} - \ov{z} \, , \, d \zeta \right>}{|\zeta -z |^{2}} 
   \we \frac{\left< (EW)_{\ve'}, d \zeta  \right>}{\left< (EW)_{\ve'} \, , \, \zeta -z \right>}  \\ \nonumber
   &\quad \we \sum_{i+j=n-2} \left[ \frac{\left< d\ov{\zeta} - d\ov{z} \, , \, d \zeta \right>}{|\zeta -z |^{2}} \right]^{i} \we \left[  \dbar_{\zeta,z}  \frac{\left< (EW)_{\ve'}  ,  d \zeta  \right>}{\left< (EW)_{\ve'}, \zeta -z \right>} \right]^{j}. 
\end{align} 
Then the homotopy identity holds
\[
  \Om^{0}_{0,q} - \Om^{(EW)_{\ve'}}_{0,q} = \dbar_{\zeta} \Om^{0, (EW)_{\ve'}}_{0,q} + \dbar_{z} \Om^{0, (EW)_{\ve'} }_{0, q-1}, \quad \text{for $(z, \zeta) \in D' \times b D_{j}$}.  
\]
We have for $ z \in D'$, 
\begin{align} \label{ibp0Wve} 
&\var_{\ve}(z) = \dbar B_{q; \, D_{j}} (\var_{\ve})(z)+B_{q+1; \, D_{j}} (\dbar\var_{\ve})(z) 
+ \int_{bD_{j}} \Om_{0,q}^{(EW)_{\ve'}} (z, \zeta) \we \var_{\ve}
 \\ \nonumber &\qquad - \dbar_{z} \int_{b D_{j}} \Omega_{0,q-1}^{0, (EW)_{\ve'}}(z,\zeta) \wedge \var_{\ve}(\zeta)
- \int_{b D_{j}} \Omega_{0,q}^{0, (EW)_{\ve'}}(z, \zeta)\wedge\dbar_\zeta\var_{\ve}(\zeta). 
\end{align} 
As shown in (i), $B_{q; \, D_{j}} (\var_{\ve})$ and $B_{q+1; \, D_{j}} (\dbar \var_{\ve})$ converge to $B_{q} \var$ and $B_{q+1} \dbar \var$ respectively in $W^{1, \frac{2n}{2n-1}}(D')$-norm as $\ve \to 0$ and $j \to \infty$. By estimate \re{trest} and \re{varcvgDj}, 
\eq{trconv} 
  \| Tr(\var_{\ve}) - Tr(\var) \|_{L^{1}(bD_{j})} \leq C \| \var_{\ve} - \var \|_{W^{1,1}(D_{j})} \overset{\ve \to 0}{\ra} 0
\eeq
where $C$ can be chosen independent of $j$. Also $\Om^{(EW)_{\ve'}}(z, \zeta)$ converges uniformly to $\Om^{EW}(z, \zeta)$ on $D' \times bD_{j}$ as $\ve' \to 0$. Hence we have
\gan
  \int_{bD_{j}} \Om_{0,q}^{(EW)_{\ve'}} (z, \zeta) \we \var_{\ve} \overset{\ve, \ve' \to 0}{\ra} \int_{bD_{j}} \Om_{0,q}^{EW} (z, \zeta) \we Tr(\var) 
\end{gather*}
uniformly on $z \in D'$. Since $\Om^{EW}(z, \zeta)$ is uniformly bounded in the first variable and uniformly continuous in the second variable for $z \in D'$ and $\zeta$ in a small neighborhood of $bD$, applying \rl{trbdycvg} we get
\[ 
  \int_{bD_{j}} \Om_{0,q}^{EW} (z, \zeta) \we Tr(\var) \overset{j \to 0}{\ra}  \int_{bD} \Om_{0,q}^{W} (z, \zeta) \we Tr(\var)
  = 0, \quad (q \geq 1).  
\] 
where the convergence is uniform for $z \in D' $. This shows that the third term in \re{ibp0Wve} converges to $0$ as $\ve, \ve' \to 0$ and $j \to \infty$. Similarly, by taking the limit as $\ve, \ve' \to 0$ and then $j \to \infty$, we can show
\[
  \int_{b D_{j}} \Omega_{0,q-1}^{0, (EW)_{\ve'}}(z,\zeta) \wedge \var_{\ve}(\zeta) \ra 
  \int_{b D} \Omega_{0,q-1}^{0, EW}(z,\zeta) \wedge Tr(\var)(\zeta), 
\] 
\[ 
  \int_{b D_{j}} \Omega_{0,q}^{0, (EW)_{\ve'}}(z, \zeta)\wedge\dbar_\zeta\var_{\ve}(\zeta) \ra
  \int_{b D} \Omega_{0,q}^{0, EW}(z, \zeta)\wedge Tr(\dbar_\zeta \var) (\zeta), 
\]
where the convergence is uniform on $z \in D'$. 
Putting together above results we obtain 
\begin{align*} 
\var(z)&= \dbar B_{q} (\var)(z)+B_{q+1} (\dbar \var)(z) 
\\ &\quad
- \dbar_{z} \int_{b D} \Omega_{0,q-1}^{0, W}(z,\zeta)\wedge Tr(\var) 
- \int_{b D}\Omega_{0,q}^{0, W}(z, \zeta)\wedge Tr(\dbar_\zeta \var), \quad z\in D'
\end{align*} 
in the distribution sense. 

Finally we check \re{hfibp1} and \re{hfibp2}. Write $ D_{c} = U \sm D$. 
Let $(EW)_{\ve'}$ and $\Om^{0, (EW)_{\e}}$ be defined as in \re{EWapprox} and \re{Omapprox}. Set $\phi = E \var$ or $E \dbar \var$, so $\phi \in W^{1, 1}(U)$. By \rp{stokesdist},  we have for $z \in D'$, 
\begin{align} \label{hfibpC1W}
&\int_{b D_c}  \Om^{0, (EW)_{\ve'}} (z, \zeta) \we Tr(\phi)
= \int_{D_c}  d \left( \Om^{0, (EW)_{\ve'}} (z, \zeta) \we \phi \right) \\ \nonumber
& \qquad = \int_{D_c} \dbar_{\zeta}  \Om^{0, (EW)_{\ve'}} \we \phi + \int_{D_c} \Om^{0, (EW)_{\ve'}} \we \dbar_{\zeta}  \phi
\\ \nonumber & \qquad = \int_{D_c} \Om^{0} (z, \zeta) \we \phi  -  \int_{D_c} \Om^{(EW)_{\ve'}} (z, \zeta) \we \phi
\\ \nonumber &\qquad \qquad + \, \dbar_{z} \int_{D_{c}} \Om^{0, (EW)_{\ve'}} (z, \zeta) \we \phi  + \int_{D_{c}} \Om^{0, (EW)_{\ve'}} (z, \zeta) \we \dbar \phi, 
\end{align}
 As $\ve' \to 0$, the $\Om^{(EW)_{\ve'}}$, $\Om^{0, (EW)_{\ve'}}$ converge uniformly to $ \Om^{EW} = \Om^{W}$ and $\Om^{0, EW} = \Om^{0, W}$ for $(z, \zeta) \in D' \times \ov{D_{c}}$, respectively. Thus 
\begin{align*}
  \int_{D_{c}} \Om^{(EW)_{\ve'}} (z, \zeta) \we \phi 
  &\overset{\ve' \to 0}{\ra} \int_{D_{c}} \Om^{W} (z, \zeta) \we \phi =0. 
\end{align*} 
Letting $\ve' \to 0$ in \re{hfibpC1W}  we get
\begin{align} \label{velimit0}
\int_{b D_{c}}  \Om^{0, W} (z, \zeta) \we  \phi
& = \int_{D_{c}} \Om^{0} (z, \zeta) \we \phi  - \dbar_{z} \int_{D_{c}}  \Om^{0, W} (z, \zeta) \we \phi   \\ \nonumber
&\quad+ \int_{D_{c}} \Om^{0, W} (z, \zeta) \we \dbar \phi. 
\end{align}		
in the distribution sense. This completes the proof of formula \re{homf1} for $bD \in C^1$, $W \in C^{1} (D \times (U \sm D))$ and $\var, \dbar \var \in W^{1,1} (D)$.  
\end{proof}

The key to our estimate is the control of the blow-up order of derivatives of the Leray map $W(z, \zeta)$ as $\zeta$ approaches the boundary from outside the domain. 
Let $D$ be a bounded domain in $\C^{n}$. Define for $\del > 0$, 
\[ 
D_{\del} = \{ z \in \C^{n}: \; \dist (z, \ov{D}) < \del  \}, \quad D_{- \del} = \{z \in D: \dist(z, b D) > \del \}. 
\] 
Gong \cite{GX18} proved the following result: 
\pr{regLeray}
Let $D$ be a bounded domain in $\C^{n}$ with $C^{2}$ boundary. Let $\rho_{0}$ be a $C^{2}$ defining function of $D$. That is, there exist a neighborhood $\mc{U}$ of $\ov{D}$ such that $D = \{ z \in \mc{U}: \rho_{0} < 0 \}$ and $\na \rho_0 \neq 0$ on $bD$. Then there exist a real function $\wti{\rho_{0}} \in C^{2} (\C^{n}) \cap C^{\infty} (\C^{n} \sm \ov{D})$ such that $\wti{\rho_0} = \rho_{0} $ in $\ov{D}$, and for $0 < d(x):= dist(x, D) < 1$, we have
\eq{regdfest}
  |\pa_{x}^{i} \wti{\rho_{0}} (x) | \leq C_{i} | \rho_0 |_{C^{2}(\ov{D})} (1+ d(x)^{2-i} )
\eeq
for $i = 0,1,2,\dots $. We call $\wti{\rho_{0}}$ the regularized defining function with respect to $\rho_0$. 

If in addition $D$ is strictly pseudoconvex. Let $\rho_1 = e^{L_{0} \rho_{0} -1}$, where $L_0$ is sufficiently large so that $\rho_1$ is strictly plurisubharmonic in a neighborhood $\om$ of $bD$. Let $\rho$ be the regularized defining function with respect to $\rho_1$. Then there exist $\del > 0$ and function $W$ (called regularized Leray map) in $D_{\del} \times (D_{\del} \sm D_{-\del})$ satisfying the following. 
	\begin{itemize}[leftmargin=1cm]
		\item[(i)] $W: D_{\del} \times (D_{\del} \sm D_{-\del}) \to \C^{n}$ is a $C^{1}$ mapping, $W(z, \zeta)$ is holomorphic in $z \in D_{\del}$, and $\Phi (z, \zeta) = W(z, \zeta) \cdot (\zeta -z) \neq 0$ for $\rho(z) < \rho(\zeta)$. 
		\item[(ii)] If $|\zeta -z| < \ve$, and $\zeta \in D_{\del} \sm D_{-\del} $, then $\Phi(z, \zeta) = F(z, \zeta) M(z, \zeta)$, $M(z, \zeta) \neq 0$ and 
		\[ 
		F(z, \zeta) = - \sum \DD{\rho}{\zeta_{j}} (z_{j} - \zeta_{j}) + \sum a_{jk} (\zeta) (z_{j} - \zeta_{j}) (z_{k} - \zeta_{k}), 
		\] 
		\[ 
		Re F(z, \zeta) \geq \rho(\zeta) - \rho(z) + |\zeta -z|^{2} / C, 
		\]
		with $M,F \in C^{1}(D_{\del} \times (D_{\del} \sm D_{- \del}))$ and $a_{jk} \in C^{\infty} (\C^{n})$. 
		\item[(iii)] For each $z \in D_{\del}$, $\zeta \in D_{\del} \sm \ov{D}, 0 \leq i, j \leq \infty$, the following holds:
		\eq{West}
		| \pa_{z}^{i} \pa_{\zeta}^{j} W(z, \zeta) | \leq C_{i,j}(D, |\rho_{0}|_{\ov{D}, 2}, \del) (1+ dist ^{1-j} (\zeta , D)). 
		\eeq
	\end{itemize} 
	\epr
	The corresponding holomorphic support function $\Phi (z, \zeta)= W(z, \zeta) \cdot (z- \zeta)$ satisfies the following estimate: 
	near every $\zeta^{\ast} \in b D$, there exist a neighborhood $V$ of $\zeta^{\ast}$ such that for all $z \in V$, there exist a coordinate map $\phi_{z}: V \to \R^{2n}$ given by $\phi_{z}: \zeta \in V \to (s, t) = (s_{1}, s_{2}, t_{3},\dots, t_{2n})$. Furthermore, for $z \in V \cap D$, $\zeta \in V \sm D$: 
	\eq{Phie1}
	|\Phi (z, \zeta) | \geq c \left( d(z) + s_{1} + |s_{2}| + |t|^{2} \right), 
	\eeq
	\eq{Phie2}
	|\Phi(z, \zeta)| \geq c |z - \zeta|^{2}, \quad \quad |\zeta -z| \geq c | (s_{2}, t) |, 
	\eeq
	where $c > 0$ is a constant.  In particular, 
	\[
	  \Phi(z, \zeta) \neq 0, \quad \text{for $z \in D$ and $\zeta \in D_{\del} \sm D$. }
	\] 
\le{cdestle}
Let $D$ be a bounded domain in $\C^{n}$ with $C^2$ boundary. Let $\rho$ be the regularized defining function as in \rp{regLeray}. Assume $\pa_{\zeta_1} \rho(\zeta^{\ast})  \neq 0$, for some $\zeta^{\ast} \in bD$. Then $\phi_{\zeta^{\ast}} = (\phi^1, \phi^2, \dots, \phi^{2n})$ given by  
\begin{align}  \label{cdphidef} 
s_{1} &= \phi^{1} (\zeta) =  \rho (\zeta), \quad s_{2} = \phi^{2} (\zeta) = Im (\rho_{\zeta} \cdot (\zeta - \zeta^{\ast})), 
\\ \nonumber
\left( t_{2k-1}, t_{2k} \right) &= \left( \phi^{2k-1}(z, \zeta), \phi^{2k}(z, \zeta) \right) = \left( Re(\zeta_{k} -\zeta^{\ast}_{k}) , \, Im(\zeta_{k} - \zeta^{\ast}_{k} ) \right), \quad k = 2, \dots ,n. 
\end{align}
defines a $C^{1}$ coordinate transformation in some neighborhood $V_0$ of $\zeta^{\ast}$. Furthermore, $\phi^{-1}$ satisfies for $m = 0, 1, 2, \dots$, 
\eq{phiinveqn}
   \left|\pa_{s}^{m} \phi^{-1} (s) \right| \leq C(1+ d(\phi^{-1}(s))^{1-m} ), \quad s \in \phi(V_0 \sm D ) . 
\eeq
\ele
\begin{proof}
	From \rp{regLeray}, we have $\rho \in C^{2}(\C^{n})$, $\phi \in C^{1} (\C^{n})$. Up to a nonzero scalar multiple, the Jacobian matrix at $\zeta = \zeta^{\ast}$ is: 
	\gan
	D \phi \left. \right|_{\zeta = \zeta^{\ast}} =
	\begin{pmatrix}
		\DD{\rho}{\zeta_{1}} &  \DD{\rho}{\ov{\zeta_{1}}}&  \DD{\rho}{\zeta_{2}} &  \DD{\rho}{\ov{\zeta_{2}}} &  \DD{\rho}{\zeta_{3}} &  \DD{\rho}{\ov{\zeta_{3}}} & \hdots &  \DD{\rho}{\zeta_{n}} &  \DD{\rho}{\ov{\zeta_{n}}} \\[0.3cm]
		\DD{\rho}{\zeta_{1}} & - \frac{\pa \rho}{\pa \ov{\zeta_{1}}} & \DD{\rho}{\zeta_{2}} & -\frac{\pa \rho}{\pa \ov{\zeta_{2}}} & \DD{\rho}{\zeta_{3}} & - \DD{\rho}{\ov{\zeta_{3}}} &\cdots & \DD{\rho}{\zeta_{n}} & -\DD{\rho}{\ov{\zeta_{n}}} \\[0.3cm]
		0 & 0 & 1 & 1 & 0 & 0 & \cdots & & 0  \\[0.3cm]
		0 & 0 & 1 & -1 & 0 & 0 & 0 & \cdots & 0 \\[0.3cm]  
		\vdots & \vdots &\vdots &\vdots &\vdots &\vdots &\ddots &\ddots &\vdots \\[0.3cm]
		0 & 0 & 0 & 0 & 0 & 0 & \cdots & 1 & 1 \\[0.3cm] 
		0 & 0 & 0 & 0 & 0 & 0 & \cdots & 1 & -1 
	\end{pmatrix}_{2n \times 2n }. 
\end{gather*}
If $n=1$, $Det (D \phi) \left. \right|_{\zeta = \zeta^{\ast}}  = -2 \DD{\rho}{\zeta_{1}} \DD{\rho}{\ov{\zeta_{1}}} \neq 0$. Suppose we have proved for $k \geq 1$, denote by $D_{k} $ and $D_{k+1}$ the determinants of $D \phi \left. \right|_{\zeta= \zeta^{\ast}}$ when $n =k$ and $n= k+1$. Computing the determinant using row expansion of second to the last row $ (0, 0, \cdots, 1, 1)$ in the above matrix, we get
\[ 
D_{k+1} = - D_{k} - D_{k} = -2 D_k \neq 0. 
\]
Thus $Det (D \phi) \left. \right|_{\zeta = \zeta^{\ast}} \neq 0$ for all $n \geq 1$. By the inverse function theorem, there exist a neighborhood $V_0$ of $\zeta^{\ast}$ such that $\phi : V_0 \to \phi(V_0)$ is  a $C^{1}$ diffeomorphism and $\phi^{-1} \in C^{1} (\phi (V_0))$.

Next, we analyze the inverse of $D \phi$. Replacing the second row in the above matrix by 
\gan
\begin{pmatrix}
	\DD{\phi^2}{\zeta_{1}} & \DD{\phi^2}{\ov{\zeta_{1}}}  & \DD{\phi^2}{\zeta_{2}} & \DD{\phi^2}{\ov{\zeta_{2}}} & \DD{\phi^2}{\zeta_{3}} & \DD{\phi^2}{\ov{\zeta_{3}}} &\cdots & \DD{\phi^2}{\zeta_{n}} & \DD{\phi^2}{\ov{\zeta_{n}}}
\end{pmatrix}, 
\end{gather*}
we obtain the Jacobian matrix $D \phi$. Leaving out the constant $\frac{1}{2i}$, we compute for $i = 1,\dots,n$, 
\begin{align} \label{phi2der1}
&\DD{\phi^2}{\zeta_i} = \DD{}{\zeta_i} \left( \rho_{\zeta} \cdot (\zeta - \zeta^{\ast}) - \rho_{\ov{\zeta}} \cdot (\ov{\zeta -\zeta^{\ast}}) \right)
\\ \nonumber &\quad = \sum_j \frac{\pa^{2} \rho}{ \pa \zeta_{i} \pa \zeta_j} (\zeta) (\zeta_j - \zeta^{\ast}_j) +  \DD{\rho}{\zeta_i} (\zeta) - \sum_j \frac{\pa^{2} \rho}{ \pa \zeta_{i} \pa \ov{\zeta_j}} (\zeta) (\ov{\zeta_j - \zeta^{\ast}_j}). 
\end{align} 
\begin{align} \label{phi2der2}
&\DD{\phi^2}{\ov{\zeta_i}} = \DD{}{\ov{\zeta_i}} \left( \rho_{\zeta} \cdot (\zeta - \zeta^{\ast}) - \rho_{\ov{\zeta}} \cdot (\ov{\zeta -\zeta^{\ast}}) \right) 
\\ \nonumber &\quad = \sum_j \frac{\pa^{2} \rho}{ \pa \ov{\zeta_{i}} \pa \zeta_j} (\zeta) (\zeta_j - \zeta^{\ast}_j) - 
\sum_j \frac{\pa^{2} \rho}{ \pa \ov{\zeta_{i}} \pa \ov{\zeta_j}} (\zeta) (\ov{\zeta_j - \zeta^{\ast}_j}) - \DD{\rho}{\ov{\zeta_i}} (\zeta). 
\end{align} 
By the inverse function theorem $[D \phi^{-1}]  = [D \phi]^{-1} \circ \phi^{-1}$ in $\phi(V_0)$. Recall the formula
\eq{invformula}
A^{-1} = \frac{1}{det (A)} Adj (A), 
\eeq
where $Adj (A)$ is the adjugate of $A$. 
Set $A = D \phi$. Then the entries of $Adj(A)$ and $det (A)$ are linear combinations with constant coefficients of 
\eq{adjdetA}
\DD{\rho}{\zeta_i}  \DD{\phi^2}{\zeta_j}, \quad \DD{\rho}{\zeta_i}  \DD{\phi^2}{\ov{\zeta_j}}, \quad \DD{\rho}{\ov{\zeta_i}}  \DD{\phi^2}{\ov{\zeta_j }}, 
\eeq
where $\DD{\phi^2}{\zeta_j}$ and $\DD{\phi^2}{\ov{\zeta_j}}$ are given by \re{phi2der1} and \re{phi2der2}.  
In view of \re{phi2der1} and \re{phi2der2}, these expressions are products of the form 
\eq{adjdetA2}
(D\rho)(\zeta) (D^{2} \rho (\zeta)) N(\zeta - \zeta^{\ast}), \quad (D\rho)(\zeta) (D \rho (\zeta)) N(\zeta - \zeta^{\ast}), 
\eeq
where $D \rho$ and $D^{2} \rho$ denote the first and second derivatives of $\rho$ and $N(\zeta - \zeta^{\ast})$ takes the form $\zeta_{j} - \zeta^{\ast}_{j}$ or $\ov{\zeta_{j} -\zeta^{\ast}_{j}}$. 
By \re{invformula} the entries of $[D \phi^{-1}] (s) = [D \phi]^{-1} \circ \phi^{-1}(s) $ take the form $  \frac{P(\zeta)}{Q(\zeta)} \circ \phi^{-1} $, where $Q(\zeta) \circ \phi^{-1} \neq 0$ in $\phi(V_0)$, and $P(\zeta)$ and $Q(\zeta)$ are some linear combination of expressions in \re{adjdetA2}.  

By \re{regdfest} the following estimates hold for $ \zeta \in V_0 \cap (\C^{n} \sm D)$: 
\eq{rhophiest}
\left| \pa_{\zeta}^{i} \rho(\zeta) \right| \leq C \left( 1 + d(\zeta)^{2-i} \right), \quad \left| \pa_{\zeta}^{i} \phi(\zeta) \right| \leq C \left( 1 + d(\zeta)^{1-i} \right), 
\eeq 
where $d(\zeta) = dist(\zeta, D)$.  We show that $\phi^{-1}$ satisfies the estimate :
\eq{phiinvest}
\left|\pa_{s}^{m} \phi^{-1} (s) \right| \leq C(1+ d(\phi^{-1}(s))^{1-m} )
\eeq
for $s \in \phi(V_0 \cap (\C^{n} \sm D))$ and $m = 0,1,2, \dots$.  Since $\phi^{-1} \in C^{1}(\phi(V_0))$, \re{phiinvest} holds for $m =1$. We have
\[ 
\pa_s \phi^{-1} (s) = \frac{P(\zeta)}{Q(\zeta)} \circ \phi^{-1} (s). 
\]
Applying chain rule we get, 
\begin{align*}
\pa^{2}_s \phi^{-1}(s) &= \frac{ (\pa_{\zeta} P)Q  - P (\pa_{\zeta} Q) }{Q^{2}} (\phi^{-1}(s)) \cdot \pa_s \phi^{-1} (s) \\
&= \frac{ (\pa_{\zeta} P)QP - P^{2}(\pa_{\zeta} Q) }{Q^{3}} (\phi^{-1} (s)). 
\end{align*}
In general, we can write $\pa_{s}^{m} \phi^{-1} (s) $ as a finite linear combination of 
\eq{PQzetader}
\frac{ \left[ \pa_{\zeta}^{j_{1}} P \cdots \pa_{\zeta}^{j_l} P  \right] \left[ \pa_{\zeta}^{k_{1}} Q \cdots \pa_{\zeta}^{k_{l'}} Q \right] P^{m_{1}} }{Q^{m_{2}}} \circ \phi^{-1} (s) 
\eeq
\[ 
   \sum_{l} j_{l} + \sum_{l'} k_{l'} = m -1, \quad m_{1}, m_{2} \geq 1. 
\]
Since $P$ and $Q$ are linear combinations of expressions in \re{adjdetA2}, by the first inequality in \re{rhophiest} we obtain that the expression in \re{PQzetader} is bounded by
\[ 
C (1 + d(\phi^{-1} (s))^{2 - 2 - (m-1)} ) = C (1 + d(\phi^{-1} (s))^{1-m}), 
\]
for $s \in \phi(V_0 \cap (\C^{n} \sm D))$ and $m = 1, 2, \dots $. This proves \re{phiinveqn}.  
\end{proof}

We now construct the coordinate system $(V, \phi)$ mentioned in the remark after \rp{regLeray}. Since $d \rho \neq 0$ on $bD$, by a linear change of coordinates we can assume that $\pa_{\zeta_{1}} \rho \neq 0$ at $\zeta^{\ast} \in bD$. By \rl{cdphidef} we can define $C^1$ coordinate transformation $\phi_{\zeta^{\ast}}$ (\re{cdphidef}) in some ball $B_{\ve} (\zeta^{\ast} )$ of small radius $\ve> 0$. 
We can find $\ve_0 > 0$ sufficiently small, such that $\phi_{z}$ defined by replacing $\zeta^{\ast} $ by $z$ in \re{cdphidef}  is  a $C^{1}$ coordinate transformation in $B_{\ve_0} (z)$, for all $z$ in some neighborhood $\om_{\zeta^{\ast} }$ of $\zeta^{\ast} $. Define $V = \om_{\zeta^{\ast} } \cap B_{\ve_{0}/2} (\zeta^{\ast} )$, then $| \zeta -z | < \ve_0$ for $z, \zeta \in V$, and thus $\phi_z$ defines a coordinate transformation on $V \subset B_{\ve_0} (z)$. 

We end the section with a trivial estimate for the top form $q =n$:  
\pr{topfest}
  Let $D$ be a bounded domain in $\C^{n}$ whose boundary is locally given by graphs of Liptschitz functions. Let $k \geq 0$ be an integer. Suppose that $\var$ is a $\dbar$-closed $(0,n)$- form. Then there exist a linear operator $S$ so that $\dbar S \var = \var$ and 
\eq{topfesteqn}
  \| S \var \|_{W^{k+1,p}(D)} \leq C(n,p) \| \var \|_{W^{k,p}(D)}. 
\eeq
\epr
\begin{proof}
   Let $B_R(0)$ be some ball centered at $0$ of radius $R$ such that $D \subset \subset B_R(0)$. Extend each component of $\var$ to a $W^{k,p}(B_R(0))$ function with compact support in $B_R(0)$.  Denote the resulting extended form by $\wti{\var}$. Since $\wti{\var}$ is a $(0,n)$-form, $\wti{\var}$ is $\dbar$-closed. Applying the homotopy formula for $B_{R} (0)$ (see \cite{WS89}*{p.~314}) and \rp{u0est}, we obtain the desired estimate. 
\end{proof}
 
\section{Estimates for $H_q$}  
We first prove a lemma which will be used in our main estimate. 
\le{intestl}
Let $0 < \del < \yh$. \\
(i) We have
\eq{intest0}
  \int_{0}^{1}  \int_{0}^{1}  \frac{s^{1+ \all} \, dt \, ds}{(\del + s +t^{2} )^{3}} \leq
   \begin{cases}
  	C(\all) \del^{\all - \yh}  &  \text{if $ 0 \leq \all < \yh$},  \\ 
  	C(\all ) ( 1+ | \log \del |) &  \text{if $\all  = \yh$},  \\
  	C(\all) &  \text{if $\all > \yh$}. 
  \end{cases} 
\eeq
(ii) If $0 < \all < 1$, we have
\eq{intest1}
 \int_{0}^{1} \int_{0}^{1}  \int_{0}^{1} \frac{s_{1}^{\all -1} t^{2n-3} \,  ds_{1} \, ds_{2} \, dt  }{(\del + s_{1} + s_{2} + t^{2})^{3}(\del + s_{1} + s_{2} + t)^{2n - 3}} \leq C (n, \all) \del^{- \frac{3}{2} +\all }. 
\eeq
\ele 
\begin{proof}
(i) Denote the integral by $I$ and split the domain of integration into three regions. \\
$R_{1}: \del + s > t > t^{2}$. We have
\[
  I \leq \int_{s=0}^{1}  \int_{t=0}^{\del + s} \frac{s^{1+ \all} \, dt \, ds}{(\del + s )^{3}}  
  \leq  \int_{s=0}^{1} (\del + s)^{-1 + \all} \, ds \leq \frac{C}{\all}. 
\]
$R_{2}: t^{2} < \del + s < t$. We have
\begin{align*}  
  I \leq \int_{s=0}^{1} \int_{t= 0}^{\sqrt{\del + s}} \frac{s^{1+ \all} \, dt \, ds}{(\del + s)^{3}} \leq \int_{0}^{1} (\del +s)^{-\frac{3}{2} + \all} \, ds  \\ 
  \qquad \leq 
 \begin{cases}
  C(\all) \del^{\all - \yh}  &  \text{if $ 0 \leq \all < \yh$},  \\ 
  C(\all ) (1 + \left| \log \del \right|) &  \text{if $\all  = \yh$},  \\
  C(\all) &  \text{if $\all > \yh$}. 
  \end{cases} 
\end{align*}
$R_{3}: \del + s < t^{2} < t$. We have
\begin{align*}  
I \leq \int_{s=0}^{1} \int_{t= \sqrt{\del + s}}^{1} \frac{s^{1+ \all} \, dt \, ds}{t^{6}} \leq \int_{0}^{1} (\del +s)^{-\frac{3}{2} + \all} \, ds  \\ 
\qquad \leq 
\begin{cases}
C(\all) \del^{\all - \yh}  &  \text{if $ 0 \leq \all < \yh$},  \\ 
C(\all ) ( 1 + |\log \del |) &  \text{if $\all  = \yh$},  \\
C(\all) &  \text{if $\all > \yh$}. 
\end{cases} 
\end{align*} 
Put together the estimates we obtain \re{intest0}. \\ 
(ii) Denote the integral by $I$, and split the domain of integration into seven regions. \\
$R_{1}: t> t^{2} > \del, s_{1}, s_{2}$. We have
\[
  I \leq \int_{\sqrt{\del}}^{1} \frac{t^{2n-3}}{t^{6}t^{2n-3}} \left( \int_{0}^{t^{2}} s_{1}^{\all -1} \, ds_{1} \right) \left( \int_{0}^{t^{2}} \, ds_{2} \right) \, dt 
  \leq C(\all) \int_{\sqrt{\del}}^{1} t^{-4 + 2 \all} \, dt \leq C(\all) \del^{-\frac{3}{2} +\all}. 
\] 
$R_{2}: t> \del > t^{2}, s_{1}, s_{2}$. We have
\[ 
   I \leq  \del^{-3} \left( \int_{\del}^{\sqrt{\del}} \frac{t^{2n-3} }{t^{2n-3}} \, dt \right) \left( \int_{0}^{\del} s_{1}^{\all-1} \, ds_{1} \right) \left( \int_{0}^{\del} \, ds_{2} \right)
   \leq C(\all) \del^{-\frac{3}{2} + \all}. 
\] 
$R_{3}: t> s_{1} > \del, t^{2}, s_{2}$. We have
\[
   I \leq \int_{\del}^{1} \frac{s_{1}^{\all - 1}}{s_{1}^{3}} \left( \int_{0}^{\sqrt{s_{1}}} \frac{t^{2n-3}}{t^{2n-3}} \, dt \right) \left( \int_{0}^{s_{1}} \, ds_{2} \right) \, ds_{1}  \leq C \del^{-\frac{3}{2} + \all}. 
\]
$R_{4}: t> s_{2} > \del, t^{2}, s_{1}$. We have
\[ 
  I \leq \int_{\del}^{1} \frac{1}{s_{2}^{3} } \left( \int_{0}^{\sqrt{s_{2}}} \frac{t^{2n-3}}{t^{2n-3}} \, dt \right) 
  \left( \int_{0}^{s_{2}} s_{1}^{\all -1} \, ds_{1} \right) \, ds_{2}   \leq C(\all) \del^{-\frac{3}{2} + \all}.  
\]
$R_{5}: \del > t, t^{2}, s_{1}, s_{2}$. We have
\[ 
   I \leq \del^{-3} \del^{-(2n-3)} \left( \int_{0}^{\del} t^{2n-3} \, dt \right) \left( \int_{0}^{\del} s_{1}^{\all -1} \, ds_{1} \right)  \left( \int_{0}^{\del} \, ds_{2} \right) 
   \leq C(n,\all) \del^{-1 + \all}. 
\]
$R_{6}: s_{1} > \del, t, t^{2}, s_{2}$. We have
\[ 
   I \leq \int_{\del}^{1} \frac{s_{1}^{\all -1}}{s_{1}^{3} s_{1}^{2n-3}} \left( \int_{0}^{s_{1}} t^{2n-3} \, dt \right) \left( \int_{0}^{s_{1}} \, ds_{2} \right) \, ds_{1} 
   \leq C(n,\all) \del^{-1 + \all}. 
\]
$R_{7}: s_{2} > \del, t, t^{2}, s_{1}$. We have
\[
  I \leq \int_{\del}^{1} s_{2}^{-3} s_{2}^{-(2n-3)} \left( \int_{0}^{s_{2}} t^{2n-3} \, dt \right) \left( \int_{0}^{s_{2}} s_{1}^{\all -1} \, ds_{1} \right)  \, ds_{2} 
  \leq C(n,\all)\del^{-1 + \all}. 
\]
Put together the estimates we obtain \re{intest1}. 

For $q \geq 1$, we can write the solution operator as
\eq{solopt} 
H_{q} \var  = u_{0} + u_{1}, 
\eeq
where 
\eq{u0u1}
   u_{0} (z) = \int_{U} \Om^{0}_{0,q-1} (z,\zeta) \we E \var , \quad u_{1} (z)= \int_{U \sm D}\Om^{01}_{0, q-1} (z, \zeta) \we [\dbar, E] \var. 
\eeq
\pr{u0est}
Let $1 < p < \infty$,  and let $U$ be a domain in $\C^{n}$, with $n > 1$. Let $u_0$ be defined as in \re{u0u1}. 
Suppose $\var \in W^{k,p}(U)$, for some nonnegative integer $k$. Then $u_{0} \in W^{k+1, p} (U)$, and
\eq{U0esteqn}
\| u_{0} \|_{W^{k+1,p}(U)} \leq C(n,p) \| \var \|_{W^{k,p} (U)}. 
\eeq
\epr
\begin{proof} 
	Let $f$ be a coefficient function of $\var$, up to a constant, $u_0$ can be written as a finite linear combination of 
	\[  
	\int_{U} \frac{\ov{\zeta^{i} - z^{i}} }{|\zeta -z|^{2n}} f(\zeta) \, dV(\zeta) 
	= \frac{1}{n-1} \int_{U} \pa_{z_{i}} \left(| \zeta -z|^{2-2n} \right) 
	f(\zeta) \, dV(\zeta) = c_{0} \pa_{z_{i}} Nf(z), 
	\]
	where $N$ denotes the Newtonian potential. Thus we just have to show that 
	\[ 
	\| Nf \|_{W^{k+2,p}(U)} \leq C(n, p) \| f \|_{W^{k,p}(U)} . 
	\]
	The proof is by Calderón-Zygmund theory. The $k=0$ case is proved in Theorem 9.9 in \cite{G-T01}*{p.~230}. Assume $k \geq 1$, we would like to move the derivatives onto $f$. Since $f$ is compactly supported in $U$, we can trivially extend $f$ to a function $\ti{f}$ in $W^{k,p}_{0}(\C^{n})$. Denoting by $\Gm$ the kernel of the Newtonian potential, we have
	\begin{align*}
	D_{x}^{k} N f(x) 
	&=  D_{x}^{k} \int_{\C^{n}}  \Gm(x-y) \ti{f}(y) \, dy  \\
	&= D_{x}^{k} \int_{\C^{n}} \Gm(y) \ti{f}(x-y) \, dy \\
	&= (-1)^{k} \int_{\C^{n}} \Gm(y) D_{y}^{k} \ti{f}(x-y) \, dy \\
	&= (-1)^{k} \int_{U} \Gm(x-y) D_{y}^{k} f(y) \, dy 
	= (-1)^{k} N(D^{k} f). 
	\end{align*} 
	Thus 
	$\left\| D^{k+2} Nf \right\|_{L^{p}(\Om)} = \left\| D^{2} N(D^{k} f) \right\|_{L^{p}(\Om)} \leq C(n,p) \| D^{k} f \|_{L^{p}(\Om)}$. 
\end{proof} 

For our estimate of $u_1$ (\re{u0u1}), we need a lemma on integration by parts. 
\le{ibple}  
   Let $D$ be a bounded domain in $\R^{N}$ with $C^{1}$ boundary. Let $\all > 0$, and $j, i_k$ be nonnegative integers. Suppose $f \equiv 0$ on $b D$. \\
  (i) Suppose $f \in C^{j+\all}(\ov{D})$. Let $g_{1}, g_{2}$ be functions in $C^{\infty}(D)$ satisfying
  \begin{gather} \label{gbloword}
    | g_k(\zeta) | \leq C |d(\zeta)|^{-i_{k}}, \qquad | \pa_{\zeta_l} g_k(\zeta) | \leq C | d (\zeta) |^{-i_{k} - 1}, \quad k = 1,2
  \end{gather} 
  for $\zeta \in D$, $d(\zeta) = \dist(\zeta, bD)$,  $1 \leq l \leq n$, and some $i_{k} \geq 0$. Furthermore $i_{k} , j$ satisfy $j \geq i_{1} + i_{2}$. 
  Then we have \\ 
  \[
    \int_{D} f \pa_{\zeta_{l}} (g_{1} g_{2}) \, dV (\zeta) = \int_{D}  f \left( \pa_{\zeta_{l}} g_{1} \right) g_{2} \, dV (\zeta) + \int_{D}  f g_{1} \left( \pa_{\zeta_{l}} g_{2} \right) \, dV (\zeta). 
  \] 
  (ii)  Suppose $f \in W^{1,1} (D) \cap C^{j+\all}(\ov{D})$, and let $g_1$ be as in (i) satisfying the estimate \re{gbloword}, such that $j \geq i_{1} $. We have
  \[
    \int_{D} f ( \pa_{\zeta_{l}} g_1) \, dV(\zeta) = - \int_{D} (\pa_{\zeta_l} f) g_1 \, dV(\zeta). 
  \]
  (iii) Let $\rho$ be a $C^{1}$ defining function of $D$. Suppose $f \in W^{1,p} (D) \cap C^{\all} (\ov{D})$, for $p >1$. Let $\phi(\zeta) =(s_{1}, \hat{s})$, $\hat{s} = (s_2, \dots, s_{2n-2})$ be a coordinate system in a neighborhood $V$ of some $p \in bD$. i.e. $\phi: V \to \phi(V)$ is a $C^{1}$ diffeomorphism. Define $\wti{f} (s) = f(\phi^{-1} (s))$ for $s \in \phi(V)$. Suppose $g$ is a function in $C^{\infty}(\phi(D \cap V))$  satisfying
  \eq{ibplog} 
    | g (s) | \leq C \left( 1+ \left| \log s_1 \right| \right), \quad \left| \pa_{s_{1}} g(s) \right| \leq C s_{1}^{-1}.   
  \eeq
  for all $s_1 < 1$. 
  Then 
  \[ 
     \int_{\phi(D \cap V)} \wti{f} (s) \pa_{s_{1}} g(s) \, ds =  - \int_{\phi(D \cap V)} (\pa_{s_{1}} \wti{f} (s) )  g(s)  \, ds. 
  \]
\ele
\noindent $\mathit{Proof.}$
(i) Let $D_{-\del} = \{ z \in D: d(z) > \del \} $, with $d(z) = \dist(z, bD)$. Take a sequence of cut-off functions $\chi_{n} \in C^{\infty}_{0} (D_{-\frac{1}{n}})$ such that $0 \leq \chi_{n} \leq 1$, $\chi_{n} \equiv 1$ on $D_{- \frac{2}{n}}$ and $| \na \chi_{n} (\zeta) | \leq C |d(\zeta)|^{-1}$ for $\zeta \in D$.

It suffices to show that 
\[ 
   \int_{D} f \pa_{\zeta_{l}} (g_{1} g_{2}) (1- \chi_{n}) \, dV (\zeta) \overset{n \to 0}{\ra} 0;
\]
\[ 
   \int_{D} f (\pa_{\zeta_{l}} g_{1} ) g_{2} (1- \chi_{n}) \, dV (\zeta) \overset{n \to 0}{\ra} 0;
\]   
\[ 
    \int_{D} f g_{1} (\pa_{\zeta_{l}}g_{2} )  (1- \chi_{n}) \, dV (\zeta) \overset{n \to 0}{\ra} 0. 
\]	 
Since $f$ vanishes on $bD$ and $f \in C^{j+\all}(\ov{D})$, then $ |f(\zeta)| \leq C d(\zeta) ^{j+\all}$, for $\zeta \in D$. In view of \re{gbloword} and that $j > i_1 + i_2$, 
the integrands in the above expression are bounded above in absolute value by a positive constant times $| d(\zeta) |^{-1 + \all} \in L^{1} (D)$. Since $1 - \chi_{n}$ converges to $0$ pointwise on $D$, the result follows from the dominated convergence thereom.  
\\ \\
(ii) Let $\chi_n$ be defined as above. It suffices to show that 
\gan
   \int_{D} f \pa_{\zeta_{l}} ((1 - \chi_{n}) g_1) \, dV (\zeta)  \overset{n \to 0}{\ra} 0; 
\\
  \int_{D} (\pa_{\zeta_{l}} f) (1 - \chi_n) g_1 \, dV (\zeta) \overset{n \to 0}{\ra} 0. 
\end{gather*}
The first statement follows from the dominated convergence theorem applied to the estimate $|f \pa_{\zeta_{l}} ((1 - \chi_{n}) g_1 )| \leq C d(\zeta)^{j+\all - (i_1+1)} \leq C d(\zeta)^{-1 + \all} \in L^{1} (D)$. For the second statement, there are two cases. If $j  \geq 1$, we have $| (\pa_{\zeta_{l}} f) (1 - \chi_n) g_1 | \leq C d(\zeta)^{j -1 +\all - i_1} \leq C d(\zeta)^{-1 + \all} \in L^{1} (\ov{D})$. If $j= i_1 = 0$, then $| (\pa_{\zeta_{l}} f) (1 - \chi_n) g_1 | \leq C | \pa_{\zeta_{l}} f | \in L^{1}(D)$, by the assumption that $f \in W^{1,1}(D)$. 
\\ \\
(iii) We can assume that $f$ is compactly supported in $V$, and thus $\wti{f} = f(\phi^{-1} (s))$ is compactly supported in $\phi(V)$. Let $ \{ \chi_{n} \}$ be defined as above. Define $\wti{\chi_n} (s)  = \chi_n (\phi^{-1} (s))$ for $ s \in \phi(V)$. Then $1- \wti{\chi_n} \equiv 0$ on $\phi(D_{-\frac{2}{n}}) \cap \phi(V) $. Since $f \equiv 0$ on $bD$ and $f \in C^{\all}(\ov{D})$, we have $ |\wti{f}(s)| \leq C s_{1}^{\all}$ and 
\[
  \int_{\phi(D \cap V)} | \wti{f}(s) |  \left| \pa_{s_{1}} \left[ (1- \wti{\chi_{n} }) g(s) \right] \right| \, ds
  \leq C \int_{ \phi (D \sm D_{-\frac{2}{n}} ) \cap \phi(V)} s_{1}^{-1 + \all - \ve} \, ds \overset{n \to 0}{\ra} 0. 
\] 
where $\ve >0$ is some arbitrary small number. By Hölder's inequality, we have
\begin{gather*}
  \int_{\phi(D \cap V)} | \pa_{s_{1}} \wti{f} |  | (1 - \wti{ \chi_n}) |  \left| g(s_{1}) \right| \, ds
  \leq C \left[ \int_{ \phi( D \sm D_{-\frac{2}{n}} ) \cap \phi(V) } | \pa_{s_{1}} \wti{f} |^{p} \, ds \right]^{\frac{1}{p}} 
  \\ \qquad \qquad \qquad
  \left[ \int_{\phi( D \sm D_{-\frac{2}{n}} )  \cap \phi(V) } \left( (1 - \wti{\chi_n} ) ( 1 + | \log s_{1}| )  \right)^{p'} \, ds \right]^{\frac{1}{p'}}, 
\end{gather*}  
which converges to $0$ since $\wti{f} \in W^{1,p} (\phi(D) \cap \phi(V))$, for $p>1$.
\end{proof} 

We are now ready for the proof of our main theorem. 
\th{mt}
   Let $D \subset \subset \C^{n}$ be a bounded strictly pseudoconvex domain with $C^{2}$ boundary. For $q \geq 1$, let $H_q \var$ be given by \re{solopt}-\re{u0u1}. \\
   (i) Let $1 < p < \infty$. Suppose $\var \in W^{1,p}(D)$. Then $ H_q \var \in W^{1,p}_{\beta}(D)$, for any $0 < \beta < \yh$, and 
   \[ 
      \| H_q \var \|_{W^{1,p}_{\beta} (D)} \leq C(D, p, \beta) \| \var \|_{W^{1,p} (D) }. 
   \] 
   (ii) Let $k \geq 2$, and $2n < p < \infty$. Suppose $\var \in W^{k,p}(D)$. Then $H_{q} \var \in W^{k,p}_{\beta}(D)$, for any $0 < \beta < \yh$, and
   \[ 
     \| H_{q} \var \|_{W^{k,p}_{\beta} (D)} \leq C(D, p, \beta) \| \var \|_{W^{k,p} (D) }. 
   \] 
\eth
\noindent $\mathit{Proof.}$
(i) We have $H_q \var = u_0 + u_1 $, where $u_0$ and $u_1$ are given by formula \re{u0u1}. By \rp{u0est}, $u_{0} \in W^{k+1,p}(D)$, and the following estimate holds: 
\[ 
   \| u_{0} \|_{W^{k+1,p}(D)} \leq C(n,p) \| E\var \|_{W^{k,p} (U)} \leq C(n,p, D) \| \var \|_{W^{k,p} (D)}. 
\]
So we only need to estimate $u_{1}$. Choose $U = D_{\del}$ as in \rp{regLeray}. We will estimate
\eq{k=1int}
  F(z) = \int_{D} d(z)^{\gm p} \left|  D_{z}^{2} \int_{U \sm D} \Om^{01}_{0,q} (z, \zeta) \we [\dbar, E] \var (\zeta) \, dV(\zeta)   \right|^{p} \, dV(z), 
\eeq
where we set $\gm = 1 - \beta$. For $z \in D$, we estimate 
\[ 
   \left| D_{z}^{2}  \int_{U \sm D} \Om^{01}_{0,q} (z, \zeta) \we [\dbar, E] \var \, dV(\zeta)  \right| , 
\] 
where in the definition of $\Om^{0, W}$ (\re{Om01}) we set $W$ to be the regularized Leray map in \rp{regLeray}.
We can write the above integral as a linear combination of 
\gan
   Kf (z) := \int_{U \sm D} f_{1}(z, \zeta) \frac{\ti{N}_{\la} (\zeta -z)  }{\Phi^{n-l} (z, \zeta)} \, dV(\zeta), \quad \quad 1 \leq l \leq n-1, \quad  
\end{gather*}
where
\eq{f1W1}
  f_{1}(z, \zeta) = f(\zeta) P_{1}(W_{1}(z,\zeta), z, \zeta),  \quad W_{1} = (\hat{D}_{\zeta}W, \pa_{z}^{p_{0}}  \hat{D}_{\zeta}W(z,\zeta)), \quad p_{0} \leq 2,
\eeq
\eq{Phi}
   \Phi(z, \zeta) = W(z, \zeta) \cdot (\zeta -z), \quad \ti{N}_{1- 2l} (\zeta - z) = \frac{N_{1} (\zeta - z)}{|\zeta - z|^{2l}}. 
\eeq
Here $f$ is a coefficient function of $[\dbar, E] \var$, and $f \equiv 0 $ on $\ov{D}$. $P_{1}(w)$ denotes a polynomial in $w$ and $\ov{w}$, $\hat{D}_{\zeta}W$ denotes $W$ and its first-order $\zeta$ derivatives, and $N_{i}$ denotes a monomial of degree $i$ in $\zeta - z$ and $\ov{\zeta -z}$. $N_{i}$ and $P_{1}$ may differ when they recur. 

Let $V$ be a small neighborhood of a fixed boundary point $\zeta^{\ast} \in bD$, as given in the remarks after \rp{regLeray}. By a linear change of coordinates we can assume that $\pa_{\zeta_{1}} \rho (\zeta^{\ast}) \neq 0$. For $z \in V$, let $\phi_{z}: V \to \phi(V)$ be the coordinate transformation given by \re{cdphidef}. 
Using a partition of unity in $\zeta$ space and replacing $f$ by $\chi f$ for a $C^{\infty}$ cut-off function $\chi$, we may assume
\[
  \text{supp}_{\zeta} \, f \subset V \sm D. 
\]
Similarly by a partition of unity in $z$ space and replacing $\Om^{01}_{0,q}$ by $\chi \Om^{01}_{0,q}$ we may assume 
\[ 
   \text{supp}_{z} \, \Om^{01}_{0,q} (z, \zeta) \subset V \cap D. 
\] 
Since $\var \in W^{1,p} (D)$, we have $f \in L^{p} (U)$. By \re{West}, we have
\eq{Ntizde}
    | \pa_{z} ^{j} \ti{N}_{\la} ( \zeta -z) | \leq |\zeta - z|^{1-2l - j}, 
\eeq
\eq{Phizde}
    \pa_{z}^{j} \Phi^{-(n-l)} (z, \zeta) \leq C_{j} (D, | \rho_{0}|_{\ov{D},2}  ) | \Phi^{-(n-l) - j} (z, \zeta) | . 
\eeq
Write
\[ 
    \pa_{z}^{2} Kf (z) = \int_{U \sm D} A(z, \zeta) f(\zeta) \, dV(\zeta), 
\]
 where  $A(z, \zeta)$ is a sum of three kinds of terms:
\[ 
A_{1} (z, \zeta) = \frac{P_{1} (z, \zeta)}{\Phi^{n-l} (z, \zeta)} \pa_{z}^{2} \left\{ \ti{N}_{1- 2l} (\zeta - z) \right\}, 
\] 
\[ 
A_{2} (z, \zeta) = \frac{P_{1} (z, \zeta)}{\Phi^{n-l + 1} (z, \zeta)} \pa_{z} \left\{ \ti{N}_{1- 2l} (\zeta - z) \right\}, 
\]
\[  
A_{3} (z, \zeta) = \frac{P_{1} (z, \zeta)}{\Phi^{n-l +2} (z, \zeta)} \left\{ \ti{N}_{1- 2l} (\zeta - z) \right\}, 
\]
and $P_{1}$ has the same form as \re{f1W1}. We have
\gan
  | \pa_{z}^{2} K f (z) | \leq \int_{U \sm D} |A(z, \zeta)|^{\frac{1}{p}}  |A(z, \zeta)|^{\frac{1}{p'}} |f(\zeta)| \, dV(\zeta),
\end{gather*}
where $\frac{1}{p} + \frac{1}{p'} =1$. Apply Hölder's inequality, we get
\eq{secdrfppe}
   | \pa_{z}^{2} K f (z)|^{p} \leq \left[ \int_{U \sm D} |A(z, \zeta)| |f(\zeta)|^{p} \, dV(\zeta)  \right] \left[ \int_{U \sm D} |A(z, \zeta)| \, dV(\zeta)  \right]^{\frac{p}{p'}}. 
\eeq
By \re{Phie2} and \re{Phi}, $C' | \zeta -z | \geq | \Phi (z, \zeta) | \geq  C| \zeta -z |^{2}$. In view of \re{Ntizde} and \re{Phizde}, we have $|A_{1}| \leq C |A_{3} | $, $| A_2| \leq C | A_{3} | $, and it suffices to estimate $A_{3} (z, \zeta) $ for $l=n-1$. From now on we just take $A$ to be
\[ 
   A(z, \zeta) = \frac{P_{1}(z, \zeta)}{\Phi(z,\zeta)^{3}}  \left\{ \ti{N}_{-(2n-3)} (\zeta - z) \right\} . 
\]
By estimate \re{Phie1}, for $z \in V \cap D$ and $\zeta \in V \sm D$: 
\eq{Philbe}
    |\Phi (z, \zeta) | \geq c (d(z) + s_{1} + |s_{2} | + |t| ^{2} ), \quad |\zeta - z| \geq c | (s_{2}, t)|. 
\eeq
where $(s_{1}, s_{2}, t) = \left( \phi^{1}_{z}(\zeta), \phi^{2}_{z} (\zeta), \phi'_{z}(\zeta ) \right) $. By \re{Philbe} and integrating by polar coordinates for $s = (s_{1} = \rho, s_{2}) \in \R^{2}$ and $t = (t_{1},\dots,t_{2n-2}) \in \R^{2n-2}$, we have 
\begin{align*}
   \int_{U \sm D} |A(z, \zeta)| \, dV(\zeta)  
   &\leq C_{0}   \int_{s =0}^{1} \int_{t =0}^{1} \frac{s t^{2n-3}\, ds  \, dt}{(d (z) + s + t^{2})^{3} t^{2n-3} } \\
   &\leq C_{0} \int_{s =0}^{1} \int_{t =0}^{1} \frac{s \, ds  \, dt}{(d (z) + s + t^{2})^{3} } \\
   &\leq C_{0}' d(z)^{- \yh}  
\end{align*}
where we used \rl{intestl} (i) for the last inequality. The constant $C_{0}$ depends only on $D$, the defining function $\rho_{0}$ and is independent of $z \in D$. Using this estimate in \re{secdrfppe} we get
\begin{align} \label{Hqfdsetup}
   & \int_{D} d(z)^{\gm p} |\pa_{z}^{2} Kf(z)|^{p} \, dV(z)  
   \\ \nonumber & \qquad \leq (C_{0})^{\frac{p}{p'}} \int_{D} \int_{U \sm D} d(z)^{\gm'} |A(z, \zeta) | |f(\zeta)|^{p} \, dV(\zeta) \, dV(z)
   \\ \nonumber & \qquad \leq (C_0)^{\frac{p}{p'}} \int_{U \sm D} \left[ \int_{D}  d(z)^{\gm'} |A(z, \zeta)| \, dV(z)   \right] |f(\zeta)|^{p} \, dV(\zeta), \quad 
\end{align}
where
\eq{mtfdgm'}
  \gm' = \gm p - \left( \yh \right)  \left( \frac{p}{p'}\right) = \left( \gm - \yh \right) p + \yh. 
\eeq
For each $z \in V$, the $C^{1}$ coordinate transformation $\phi_{z}$ is given by \re{cdphidef}: 
\[ 
   \phi^{1}_{z} (\zeta) = \rho (\zeta), \quad \phi^{2}_{z} (\zeta) = Im (\rho_{\zeta} \cdot (\zeta -z)), \quad \phi'_{z} (\zeta) = \left( Re(\zeta' -z'), Im(\zeta'- z') \right). 
\] 
For $\zeta \in V$, we define $\wti{\phi}_{\zeta}: V \to \phi(V)$ to be 
\begin{gather} \label{phitizeta} 
   \wti{\phi}^{1}_{\zeta} (z) = \rho (z), \quad \wti{\phi}^{2}_{\zeta}(z) = Im (\rho_{\zeta} \cdot (\zeta -z)),  
   \\ \nonumber  \wti{\phi}'_{\zeta} (z) = \left( Re(\zeta' -z'), Im(\zeta'- z') \right), 
\end{gather}
which is a coordinate system for $z \in V$. Write $(\ti{s}_{1}, \ti{s}_{2}, \ti{t}) = (\wti{\phi}^{1}_{\zeta} (z), \wti{\phi}^{2}_{\zeta} (z), \wti{\phi}'_{\zeta} (z)) $. By \re{Philbe} we have for $z \in V \cap D$ and $\zeta \in V \sm D$, 
\begin{align} \label{Phiest}
 |\Phi (z, \zeta) | &\geq c (d(z) + \phi^{1}_{z} (\zeta) + | \phi^{2}_{z}(\zeta) | + |\phi'_{z} (\zeta)|^{2} ) \\ \nonumber
  &\geq c (d(\zeta) + |\wti{\phi}^{1}_{\zeta} (z)| + | \wti{\phi}^{2}_{\zeta}(z) | + |  \wti {\phi}' _{\zeta} (z)|^{2} ) \\ \nonumber
  &= c (d(\zeta) + |\ti{s}_{1}| + |\ti{s}_{2}| + |\ti{t} |^{2} ). 
\end{align}
and 
\eq{}
   |\zeta - z| \geq c | (s_2, t)| = c |(\wti{s}_2, \wti{t}) |.  
\eeq
Writing in polar coordinates and using that $d(z) \leq C \rho(z) = C | \ti{ s_{1}} | \leq C | \ti{s}|$ , we have by \rl{intestl} (i) again
\begin{align} \label{Hqfdest}
  \int_{D}  d(z)^{\gm'} |A(z, \zeta)| \, dV(z) 
  &\leq C \int_{\ti{s} =0}^{1} \int_{\ti{t} =0}^{1} \frac{\ti{s} ^{1 + \gm'} \wti{t}^{2n-3} \, d \ti{s}  \, d \ti{t}}{(d (\zeta) + \ti{s} + \ti{t}^{2})^{3} \wti{t}^{2n-3}} \\ \nonumber 
  &\leq 
  \begin{cases}
   C d(\zeta)^{\gm' - \yh}  &  \text{if $ 0 \leq \gm' < \yh$},  \\ 
   C ( 1 + | \log d(\zeta) |) &  \text{if $\gm' = \yh$},  \\
   C(\gm') &  \text{if $\gm' > \yh$}, 
  \end{cases}
\end{align}
If $\gm > \yh$, then by \re{mtfdgm'}, $\gm' > \yh$. Using \re{Hqfdest} in \re{Hqfdsetup}, we get 
\begin{align*}
  \left[ \int_{D} d(z)^{\gm p} |\pa_{z}^{2} Kf(z)|^{p} \, dV(z)   \right]^{\frac{1}{p}}
  &\leq C'(D, \gm') \left[ \int_{U \sm D}  |f (\zeta) |^{p} \, dV(\zeta) \right]^{\frac{1}{p}} \\
  &\leq C'(D, \gm') \| \var \|_{W^{1,p}(D)}, 
\end{align*}
Thus we have shown that
\[ 
  \| u_{1} \|_{W^{1,p}_{\beta}(D)} \leq C(D, \beta ) \| \var \| _{W^{1,p}(D)}, 
\]
for any $0 < \beta < \yh$.  
\\ \\
(ii) Next we estimate higher derivatives for $u_1$. Suppose $\var \in W^{k,p}(D)$, for $k \geq 2$ and $2n < p < \infty$. We show that $u_{1} \in W^{k,p}_{\beta}(D)$,  for any $0 < \beta < \yh$. Let $f$ be a coefficient function of $[\dbar, E] \var$. As before take $U = D_{\del}$ as in \rp{regLeray}. Then $f \in W^{k-1,p}(U)$. By \rp{Sobemb}, $f \in C^{k-2 + \all}(U)$,  $\all = 1 - \frac{2n}{p}$. Since $f \equiv 0$ on $D$, for $\zeta \in U$ the following holds: 
\eq{fvnsh}
    |\pa^{q}_{\zeta} f (\zeta) | \leq C_{q} |f|_{U; k-2+\all} \, d(\zeta)^{k-2 + \all-q}, \quad 0 \leq q \leq k-2, 
\eeq
where $d(\zeta) = \dist (\zeta, bD)$. We have
\begin{align*}
    &\int_{D} \left| \pa_{z}^{k+1} u_{1}(z) \right|^{p} d(z)^{\gm p} \, dV(z)  
    \\ \nonumber & \quad =  \int_{D} \left| \pa_{z}^{2} \int_{U \sm D} \pa_{z}^{k-1} \Om^{01}_{0,q} (z, \zeta) \we [\dbar, E] \var (\zeta) \, dV(\zeta)   \right|^{p} d(z)^{\gm p}  \, dV(z). 
\end{align*}
We can write the inner integral above as a linear combination of 
\gan
   K_{1}f(z) = \int_{U \sm D} f_{1}(z, \zeta) \frac{N_{1-\mu_{0} + \mu_{2}}(\zeta -z )}{\Phi^{n-l+\mu_{1}}(z, \zeta) |\zeta -z|^{2l+2 \mu_{2} }} \, dV(\zeta),  \quad 1 \leq l \leq n-1
   \\   f_{1}(z, \zeta) = f(\zeta) P_{1}(W_{1}(z,\zeta), z, \zeta),  \quad W_{1} = (\hat{D}_{\zeta}W, \pa_{z}^{k_{0}}  \hat{D}_{\zeta}W(z,\zeta)), 
   \\ \mu_{0} + \mu_{1} + \mu_{2} \leq k-1, \quad 1-\mu_{0} + \mu_{2} \geq 0, \quad  k_{0} \leq k-1.
\end{gather*}
We apply integration by parts in two stages. In the first stage, we integrate by parts to reduce the exponent of $\Phi$ in the denominator to $n-l$, as in Ahern-Schneider \cite{A-S79}, Lieb-Range \cite{L-R80} and Gong \cite{GX18}. See also Michel-Perotti \cite{M-P90} for estimates without using integration by parts for piecewise smooth strictly pseudoconvex domains via Seeley extension. 

Let $V$ be a small neighborhood of a fixed boundary point $\zeta^{\ast} \in bD$ as in (i). Suppose that for $z \in V \cap D$ and $\zeta \in V \sm D$, 
\[ u(z, \zeta) := \pa_{\zeta_{i_{\ast}}} \Phi(z, \zeta) \neq 0 , \quad \quad \text{for some $ i_{\ast}$. } \] 
By \re{West}, for fixed $z \in D$ the following estimates hold for $ \zeta \in V \sm D$ if $b D$ is $C^{2}$: 
\eq{Phibloword}
   |\pa_{\zeta_{i_{\ast}}}^{q} \Phi^{-k} (z, \zeta) | \leq C(D,z) \left( 1 + d(\zeta)^{1- q} \right), 
\eeq
\eq{ubloword}
  |\pa_{\zeta_{i_{\ast}}}^{q} u^{-k}(z, \zeta)| \leq C(D,z) d(\zeta)^{-q}, 
\eeq
for $q = 0, 1, 2, \dots $ and $k = 1, 2, \dots$. 
Up to a constant multiple, we rewrite $K_{1} f$ as 
\begin{align} \label{K1f}
  K_{1} f(z) &= \int_{U \sm D}  f(\zeta) h(z,\zeta) \pa_{\zeta_{i_{\ast}}} \Phi^{-(n-l+\mu_{1})}(z, \zeta) \, dV(\zeta)  \\ \nonumber
  &= \int_{U \sm D}  f(\zeta) u(z, \zeta)^{-1} h(z,\zeta) \pa_{\zeta_{i_{\ast}}} \Phi^{-(n-l+\mu_{1} -1)}(z, \zeta) \, dV(\zeta)
\end{align}
where we set
\eq{h} 
   h(z, \zeta) = P_{1} (W_{1}(z, \zeta), z, \zeta) \frac{N_{1-\mu_{0} + \mu_{2}} (\zeta  -z)}{|\zeta -z|^{2l + 2 \mu_{2}}}. 
\eeq
For fixed $z \in D$ the following holds for $\zeta \in  U \sm D$, 
\eq{hbloword}
    |\pa_{\zeta_{i_{\ast}}}^{q} h(z, \zeta)| \leq C (D,z ) d(\zeta)^{-q}, \quad q = 0, 1, 2, \dots 
\eeq
again by \re{West}. 
Then we get 
\begin{align} \label{K1fstokes}
  K_{1} f(z) &= \int_{U \sm D} f(\zeta) \pa_{\zeta_{i_{\ast}}} \left[ u^{-1}(z, \zeta) h(z, \zeta)  \Phi^{-(n-l+\mu_{1} - 1)}(z, \zeta) \right]\, dV(\zeta)
  \\ \nonumber & \quad - \int_{U \sm D} f(\zeta) \pa_{\zeta_{i_{\ast}}} \left[ u^{-1}(z, \zeta) h(z, \zeta) \right]  \Phi^{-(n-l+\mu_{1} - 1)}(z, \zeta) \, dV(\zeta)
  \\ \label{k1fstokes2} & = - \int_{U \sm D} \left[ \pa_{\zeta_{i_{\ast}}} f(\zeta)  \right] \left[ u^{-1}(z, \zeta) h(z, \zeta)  \Phi^{-(n-l+\mu_{1} - 1)}(z, \zeta) \right]\, dV(\zeta) 
  \\ \nonumber & \quad - \int_{U \sm D} f(\zeta) \pa_{\zeta_{i_{\ast}}} \left[ u^{-1}(z, \zeta) h(z, \zeta) \right]  \Phi^{-(n-l+\mu_{1} - 1)}(z, \zeta) \, dV(\zeta). 
\end{align}
We now justify the above steps. 
We have $f \in W^{1,p}(U) \cap C^{j+ \all}(U)$, for $j = k-2$.  Apply \rl{ibple} to the domain $U \sm D$ with $f \equiv 0$ on $b (U \sm D)$. By \re{ubloword}, \re{hbloword}, $u^{-1}, h$ satisfy the estimates \re{gbloword} with $i_{k} = 0$. By \re{Phie1} for fixed $z \in D$, $\left| \Phi^{- (n-l+\mu_{1} - 1)} (z, \zeta) \right| \leq C(z)$ , and 
\[ 
  \left| \pa_{\zeta_{i_{\ast}}} \Phi^{-(n-l+\mu_1 -1)} (z, \zeta) \right|  = \left| C u(z, \zeta) \Phi^{-(n-l+\mu_1)} (z, \zeta) \right| \leq C(D,z). 
\]
Thus $\Phi^{-(n-l+\mu_1 -1)}$ also satisfies the estimate \re{gbloword} with $i_{k} = 0$. Then the first equality \re{K1fstokes} follows from \rl{ibple} (i) and the second equality \re{k1fstokes2} follows from \rl{ibple} (ii).  

We can repeat this procedure $\mu_{1} (\leq k-1)$ times. Indeed, suppose we have done $m$ times, $1 \leq m \leq k-2$. Then the integral is a linear combination of terms of the form
\begin{align*}
  &\int_{U \sm D} \left( \pa_{\zeta_{i_{\ast}}}^{m_{1}} f \right) \pa_{\zeta_{i_{\ast}}}^{m_{2}} \left\{  u^{-1}, h \right\}  \Phi^{-(n-l+\mu_{1} - m )} \, dV(\zeta) \quad \quad (m_{1} + m_{2} = m)
  \\ &\qquad \qquad =  \int_{U \sm D} \left( \pa_{\zeta_{i_{\ast}}}^{m_{1}} f \right) \pa_{\zeta_{i_{\ast}}}^{m_{2}} \left\{ u^{-1}, h \right\} u^{-1} \pa_{\zeta_{i_{\ast}}}\Phi^{-(n-l+\mu_{1} - m - 1)} \, dV(\zeta), 
\end{align*} 
where $\pa^{m_2}_{\zeta_{i_{\ast}}} \left\{ u^{-1}, h \right\}$ denotes a linear combination with constant coefficients of the terms
\begin{gather} \label{u-1hdivder}
  \pa_{\zeta_{i_{\ast}}}^{\la_1} (u^{-1}) \pa_{\zeta_{i_{\ast}}}^{\la_2} (u^{-1}) \cdots \pa_{\zeta_{i_{\ast}}}^{\la_{p}} (u^{-1})
  \pa_{\zeta_{i_{\ast}}}^{\la_0} (h), 
  \\ \nonumber
  \la_{i} \geq 0 \quad \quad \sum_{i=0}^{p} \la_i = m_2 . 
\end{gather} 
Then $\pa_{\zeta_{i_{\ast}}}^{m_{1}} f \in W^{1,p}(U \sm D) \cap C^{k-2 - m_{1} +\all}(\ov{U \sm D })$. Also  $\pa_{\zeta_{i_{\ast}}}^{m_{2}} \left\{ u^{-1}, h \right\} $, $u^{-1}$ satisfy estimates \re{gbloword} for $i_{k} = m_{2}$, and $ \Phi^{-(n-l+\mu_{1} - j - 1)} \leq C(z)$. Since $k-2- m_{1} - m_{2} = k-2- m \geq 0$, the hypothesis of \rl{ibple} (i) and (ii) holds, and we can do the procedure one more time. 
 
From the above argument we can now write $K_{1} f$ (\re{K1f}) as a linear combination of 
\eq{K2f}
   K_{2}f(z) = \int_{U \sm D} \pa_{\zeta_{i_{\ast}}}^{\tau_{0}}f (\zeta) \pa_{\zeta_{i_{\ast}}}^{\tau_{1}} \left\{ u^{-1},h \right\}(z, \zeta) \Phi^{-(n-l)}(z, \zeta) \, dV(\zeta)
\eeq
\[ 
   \tau_{0} + \tau_{1} = \mu_{1},  \quad \tau_{0} < k-1, 
\]
and 
\eq{K2ti}
  K_{2}' f (z) = \int_{U \sm D} \pa_{\zeta_{i_{\ast}}}^{k-1} f(\zeta) h'(z, \zeta) u^{- (k-1)} \Phi^{-(n-l)} (z, \zeta) \, dV(\zeta), 
\eeq
where 
\[ 
   h'(z, \zeta) = P_{1} (W_{1}(z, \zeta), z, \zeta) \frac{N_{1} (\zeta  -z)}{|\zeta -z|^{2l}}, \quad \quad (\mu_{0} = \mu_{2} = 0). 
\]
In the case all $k-1$ derivatives fall onto $f$, we have the integral $K_{2}'f$. Since $\pa_{\zeta_{i_{\ast}}}^{k-1} f \in L^{p}(U \sm D)$, this reduces to the earlier $k=1$ case, and we obtain 
\eq{K2tiest}
   \int_{U \sm D} \left| \pa_{z}^{2} K_{2}' f(z, \zeta) \right|^{p} d(z)^{ \gm p} \, dV(z) \leq C \| \var \|_{W^{k,p}(D)}^{p},
\eeq
for any $\gm > \yh$. 

The above integration by parts suffices to derive the estimates in \cite{L-R80} and \cite{GX18}. For our estimates, we must go through a second stage of integration by parts for $K_2 f$ to avoid unnecessary loss in regularity. We integrate by parts with respect to the normal direction, and again we rely on the regularized Leray map. 

In view of \re{u-1hdivder} and \re{h}, we can write $  \pa_{\zeta_{i_{\ast}}}^{\tau_{1}} \left\{ u^{-1}, h \right\}(z, \zeta)$ as a linear combination of 
\begin{align} \label{divder}
  \widehat{\pa}_{\zeta_{i_{\ast}}}^{\iota} \left\{ u^{-1}, h \right\}(z, \zeta) 
  &=  \pa_{\zeta_{i_{\ast}}}^{\iota_1} (u^{-1}) \cdots \pa_{\zeta_{i_{\ast}}}^{\iota_{p}} (u^{-1})
  \\ \nonumber &\qquad \pa_{\zeta_{i_{\ast}}}^{\iota_{0}} (P_{1} (W_{1}, z, \zeta))  \pa_{\zeta_{i_{\ast}}}^{\nu_{1}} 
  \left( \frac{N_{1 - \mu_{0} + \mu_{2}}}{|\zeta -z|^{2l + 2 \mu_{2}}} \right)
\end{align} 
\eq{iotaj} 
  \sum_{j =0}^{p} \iota_j= \nu_0, \quad \nu_{0} + \nu_{1} = \tau_{1}. 
\eeq
For $z \in V \cap D$, let $\phi_{z}: U_0 \to H^{+} $ be given by \re{cdphidef}, where we denote 
\[
  H^{+} = [0,1] \times [-1,1]^{2n-1}. 
\]
For simplicity we write $\phi$ and $\phi^{-1}$ in place of $\phi_z$ and $\phi^{-1}_z$. 
Define 
\[ 
  \wti{\pa^{\tau_{0}}_{\zeta_{i_{\ast}}} f} (s) = \pa^{\tau_{0}}_{\zeta_{i_{\ast}}} f(\phi^{-1} (s)), \quad   \wti{\Phi} (z, s) = \Phi(z, \phi^{-1} (s)),  
\] 
\eq{gzs}
  g^{\iota}(z, s) = \widehat{\pa}_{\zeta_{i_{\ast}}}^{\iota} \left\{ u^{-1}, h \right\} (z, \phi^{-1} (s)) \, \left| det( D \phi^{-1}) (s) \right|, 
\eeq
where $D \phi^{-1} $ denotes the Jacobian of $\phi^{-1}$. 
Then $K_{2} f$ (\re{K2f}) can be written as a linear combination of
\begin{align} \label{K3f} 
  K_{3} f (z) 
 &= \int_{H^{+}} \wti{\pa^{\tau_{0}}_{\zeta_{i_{\ast}}} f} (s)  \frac{g^{\iota}(z, s) }{\wti{\Phi}^{n-l}(z, s) }  \, dV(s) \\ \label{ibpbefore} 
 &= \int_{H^{+}}  \wti{\pa^{\tau_{0}}_{\zeta_{i_{\ast}}} f} (s)  \pa_{s_{1}} I_{1} (z,s) \, dV(s), 
\end{align}
where
\eq{I1Hq} 
  I_{1} (z,s) = \int_{1}^{s_{1}}   
  \frac{ g^{\iota}(z, (\eta_{1}, \hat{s}))}{\wti{\Phi}^{n-l}(z, (\eta_{1}, \hat{s}))}  \, d \eta_1,
\eeq
and $  s = (s_{1}, \hat{s}), \; \hat{s} = (s_{2}, t_3 \dots, t_{2n})$. Observe that for a fixed $z \in D$, by \re{Phie1} the $\wti{\Phi}(z, s)$ and $|\phi^{-1} (s) - z|$ are bounded below by a constant depending on $d(z)$. By definition of $\phi$ and \re{regdfest}, the following holds for $\zeta \in V \sm D$: 
\eq{phibup}
   |\pa_{\zeta}^{q} \phi (\zeta) | \leq C_{q} (1 + d(\zeta)^{1-q}), \quad q = 0,1,2, \dots 
\eeq  
By estimate \re{phiinveqn}, the following holds for $s \in H^{+}$: 
\begin{align} \label{phiinvbup}
  |\pa_{s}^{q} \phi^{-1} (s) | &\leq C_{q} (1 + d(\phi^{-1} (s))^{1-q}) 
  \\ \nonumber &\leq  C_{q} (1 + s_{1}^{1-q}), \quad q = 0,1,2, \dots 
\end{align}
In particular, 
\eq{detest}
  \left| \pa_{s}^{q} [det(D \phi^{-1})] (s) \right| \leq C_{q} s_{1}^{-q}, \quad q = 0,1,2, \dots, 
\eeq
and thus 
\eq{detests1'}
   \left| \pa_{s}^{q} [det(D \phi^{-1})] ((\eta_{1}, \hat{s})) \right| \leq C_{q} \eta_{1}^{-q}, \quad q = 0,1,2, \dots  
\eeq
By \re{ubloword} and \re{hbloword}, we have 
\begin{align} \label{P1u-1est}
  &\left| \pa_{\zeta}^{q} P_{1}(z, \phi^{-1} (\eta_1, \hat{s})) + \pa_{\zeta}^{q} (u^{-1}) (z, \phi^{-1} (\eta_1, \hat{s})) \right| 
  \\ \nonumber & \qquad \leq C(D) \left[ d(\phi^{-1} (\eta_1, \hat{s} )) \right] ^{-q}  \\ \nonumber
  & \qquad \leq C(D) \left[ \phi^{1}(\phi^{-1} (\eta_1, \hat{s})) \right]^{-q} \\ \nonumber
  & \qquad = C(D) \eta_1^{-q}
\end{align}
for $q = 0,1,2 \cdots $. Applying \re{detests1'}, \re{P1u-1est} to \re{divder}, \re{gzs} we get
\begin{align} \label{giotaest}
  |g^{\iota}(z, (\eta_1, \hat{s}))| 
  &\leq C \left| \widehat{\pa}_{\zeta_{i_{\ast}}}^{\iota} \left\{ u^{-1},h \right\} (z, \phi^{-1} ((\eta_1, \hat{s}))) \right| \, \left| det( D \phi^{-1}) ((\eta_1, \hat{s})) \right| \\ \nonumber 
  &\leq C(D,z) \eta_1^{- \nu_{0}} |\phi^{-1} ((\eta_1, \hat{s})) - z|^{1- 2l - \mu_{0} - \mu_{2} - \nu_{1}}  \\ \nonumber
  &\leq C(D,z) \eta_1^{- \nu_{0}}.
\end{align} 
where $\nu_{0} = \sum_{j=0}^{p} \iota_j$. In view of \re{I1Hq} and \re{giotaest}, we have 
\begin{align*} 
| I_{1}(z, s) | \leq C(D,z) \int_{s_1}^{1} \eta_{1}^{-\nu_{0}} \, d\eta_1 
 & \leq 
\begin{cases}
 C(D, z) s_{1}^{-(\nu_{0} -1)}  & \text{if $\nu_{0} >1$. } \\ 
 C (D, z) | \log s_{1} | & \text{if $\nu_{0} =1$. } \\
 C (D,z) & \text{if $\nu_{0} = 0$. }
\end{cases}
\end{align*}
and 
\begin{align*} 
 | \pa_{s_{1}} I_{1} (z, s) | = \left| \frac{g^{\iota}  (z, (s_1, \hat{s})) }{ \wti{\Phi}^{n-l}(z,  (s_1, \hat{s}) ) } \right| 
 \leq \begin{cases}
 C(D, z) s_{1}^{- \nu_{0} }  & \text{if $\nu_{0} >1$. } \\ 
 C (D, z) s_{1}^{-1} & \text{if $\nu_{0} =1$. } \\
 C (D,z) & \text{if $\nu_{0} = 0$. }
 \end{cases}
\end{align*}
If $\nu_{0} >1$, we can apply \rl{ibple} (ii) to $K_{3}f$ (\re{ibpbefore}) with $\wti{\pa^{\tau_{0}}_{\zeta_{i_{\ast}}} f} (s) \in W^{1,p}(H^{+}) \cap C^{j+\all} (H^{+}) $ for $j = k-2 - \tau_{0} \geq 0$ ($\tau_{0} \leq k-2$) and $i_1 = \nu_{0} -1$. We have
\[ 
j-i = (k-1) - \tau_{0} - \nu_{0} \geq (k-1) - \tau_{0} - \tau_{1} \geq (k-1) - \mu_{1} \geq 0.
\]
Thus we can integrate by parts in \re{ibpbefore} to get
\eq{ibpafter}
  K_{3} f (z) = - \int_{H^{+}} \pa_{s_{1}} \wti{\pa^{\tau_{0}}_{\zeta_{i_{\ast}}} f} (s) I_{1}(z, s) \, dV(s). 
\eeq
If $\nu_0 = 1$, we can apply \rl{ibple} (iii) for $\wti{\pa^{\tau_{0}}_{\zeta_{i_{\ast}}} f} (s) \in W^{1,p}(H^{+}) \cap C^{\all} (H^{+})$,  and integrate by parts in \re{ibpbefore} to get \re{ibpafter}. 
Finally if $ \nu_0 = 0$, we can again apply \rl{ibple} (ii) with $\wti{\pa^{\tau_{0}}_{\zeta_{i_{\ast}}} f} \in C^{k-2 - \tau_0 + \all} (H^{+})$ and $i = 0$. 

We claim that we can integrate by parts in this fashion $k-1-\tau_0$ times. Suppose we did it for $m$ times, for $1 \leq m \leq (k-2) - \tau_{0} $, and we have
\gan
   K_{3} f (z) = \pm \int_{H^{+}} \pa_{s_{1}}^{m} \wti{ \pa^{\tau_{0}}_{\zeta_{i_{\ast}}}} f (s) I_{m} (z, s)  \, dV(s), 
\end{gather*}
where
\[
  I_{m} (z, s) = \underbrace{  \int_{1}^{s_{1}} \int_{1}^{\eta_1} \cdots \int_{1}^{\eta_{m-1}} }_\text{$m$ integrals}  \frac{g^{\iota} (z, ( \eta_m,  \hat{s}  ) ) \, [d \eta]^{m} }{\wti{\Phi}^{n-l}(z, (\eta_{m}, \hat{s} ) ) }, 
\]
and we denote $[d \eta]^{m} := \, d \eta_m \cdots d \eta_1$. Recall \re{phiinvbup}, 
\eq{phiinvbup1}
  |\pa_{s}^{q} \phi^{-1} (s) | \leq C_{q} (1 + s_{1}^{1-q}), \quad s \in \phi(V \sm D), \quad q = 0,1,2,\dots
\eeq
In particular $| \pa^{1}_{s} \phi^{-1} (s)| \leq C$. By the chain rule we observe that $\pa_{s_{1}}^{m} \wti{\pa^{\tau_{0}}_{\zeta_{i_{\ast}}} f} (s) = \pa_{s_{1}}^{m} \pa^{\tau_{0}} _{\zeta_{i_{\ast}}} f (\phi^{-1} (s))$ is a sum of terms of the form
\[
   \left[ \pa^{m_{0} + \tau_{0}}_{\zeta_{i_{\ast}}} f(\phi^{-1} (s)) \right]  \left[  \pa_{s}^{m_{1}} \phi^{-1} (s) \right]
   \cdots  \left[ \pa_{s}^{m_{\ell}} \phi^{-1} (s) \right] \left[ \pa_{s}^{1} \phi^{-1} \right]^{m_{\ell+1}}, 
\]
\[ 
  m_{0} + m_{1} + \cdots + m_{\ell} \leq m + 1, \quad m_{0} \leq m,  \quad m_{\ell+1} \leq m. 
\]
In view of this and \re{phiinvbup1}, we have the estimate
\begin{align} \label{fibpest}
  \left| \pa_{s_{1}}^{m} \wti{\pa^{\tau_{0}}_{\zeta_{i_{\ast}}} f} (s)\right| 
  &\leq C \left[ d(\phi^{-1} (s))  \right]^{k-2 + \all - m - \tau_{0}} \\
  &\nonumber\leq C \left[ \phi^{1} (\phi^{-1} (s)) \right]^{ k -2 + \all -m - \tau_{0}} \\
  &\nonumber \leq C s_{1}^{k - 2 + \all -m - \tau_{0}}. 
\end{align}
Write $K_{3}f$ as
\begin{align} \label{bfibp}
   K_{3} f (z) = &\pm \int_{H^{+}} \left( \pa_{s_{1}}^{m} \wti{\pa^{\tau_{0}}_{\zeta_{i_{\ast}}} f}  (s) \right) \pa_{s_{1}} I_{m+1} (z,s) \, dV(s) 
\end{align} 
where  
\gan
  I_{m+1} (z,s) =  \underbrace{  \int_{1}^{s_{1}} \int_{1}^{\eta_1} \cdots \int_{1}^{\eta_m} }_\text{$m + 1$ integrals}  \frac{g^{\iota}  (z, (\eta_{m+1}, \hat{s}) )}{\wti{\Phi}^{n-l}(z, ( \eta_{m+1}, \hat{s})) } \, [d \eta]^{m+1}. 
\end{gather*}
We have by \re{giotaest},
\begin{align*}
   \left| I_{m+1}(z,s) \right| &\leq C(D, z) \underbrace{  \int_{s_1}^{1} \int_{\eta_1}^{1} \cdots \int_{\eta_m}^{1} }_\text{$m + 1$ integrals} (\eta_{m+1} )^{-\nu_{0}} \, [d \eta]^{m+1} \\
    &\leq 
    \begin{cases}
    C(D, z) s_{1 }^{- \nu_{0} + (m+1)}, & \text{if $m+1 < \nu_{0} = \sum \iota_{j}$},   \\
    C(D,z) (1+| \log s_{1} |), & \text{if $m+1 = \nu_{0}$}, \\
    C(D,z), & \text{if $m+1 > \nu_{0}$}. 
    \end{cases}
\end{align*}
where $\iota_{j}$-s are defined in \re{iotaj}.  
\begin{align*}
   \left|  \pa_{s_{1}} I_{m+1}(z,s) \right| 
   &= \left| I_{m} (z, s) \right| \\
   &= \begin{cases}
   C(D, z) s_{1 }^{- \nu_{0} + m}, & \text{if $m+1 < \nu_{0} $},   \\
   C(D,z) s_{1}^{-1}, & \text{if $m+1 = \nu_{0}$}, \\
   C(D,z) ( 1+| \log s_{1}|), & \text{if $m+1 > \nu_{0}$}. 
   \end{cases}
\end{align*}
If $m+1 < \nu_0 = \sum \iota_{j}$, then $I_{m+1}$ satisfies estimate \re{gbloword} with $i_{k}$ replaced by $i = \nu_{0} - (m+1)$. We can apply \rl{ibple} (ii) to \re{bfibp} for $ \pa_{s_{1}}^{m} \wti{\pa^{\tau_{0}}_{\zeta_{i_{\ast}}} f} (s) \in W^{1,p}(H^{+}) \cap C^{j+\all}(H^{+})$, with $j= k-2- m - \tau_{0} \geq 0$. We have
\[ 
   j-i = (k-1) - \tau_{0} - \nu_{0} \geq (k-1) - \tau_{0} - \tau_{1} \geq (k-1) - \mu_{1} \geq 0.
\]
Thus we can integrate by parts in \re{bfibp} to get
\begin{align} \label{afibp}
      K_{3} f (z) &= \pm \int_{H^{+}}  \left( \pa_{s_{1}}^{m+1} \wti{\pa^{\tau_{0}}_{\zeta_{i}} f}  \right) (s) I_{m+1} (z, s) dV(s). 
\end{align}
If $m +1 = \nu_{0} = \sum \iota_{j}$, then $I_{m+1}$ satisfies estimate \re{ibplog}, and we can apply \rl{ibple} (iii) to obtain \re{afibp}. 
If $m + 1 > \nu_{0} = \sum \iota_{j}$, then $I_{m+1}$ satisfies estimate \re{gbloword} with $i_{k}$ replaced by $i=0$, and we again apply \rl{ibple} (ii) to obtain \re{afibp}. In conclusion, we can transform $K_{3}f$ (\re{K3f}) via integration by parts to the form 
\eq{K3f2derbefore}
  \int_{H^{+}} F (s) I_{k-1 - \tau_{0}} (z, s) \, dV(s), 
\eeq
where
$$ F(s) = \pa_{s_{1}}^{k-1 - \tau_{0}} \wti{ \pa_{\zeta_{i_{\ast}}}^{\tau_{0}}} f (s)  \in L^{p} (H^{+}), $$
and
$$  I_{k-1 - \tau_{0}} (z, s) = \underbrace{  \int_{1}^{s_{1} } \int_{1}^{\eta_1} \cdots \int_{1}^{\eta_{k-2- \tau_{0}}} }_\text{$k-1 - \tau_{0}$ integrals}  \frac{g^{\iota} (z, (\eta_{k-1 - \tau_0}, \hat{s} ) )}{\wti{\Phi}^{(n-l)}(z, (\eta_{k-1-\tau_0}, \hat{s} ) ) } \, d \eta , $$
with $d \eta = d \eta_{k-1-\tau_0} \, \cdots \, d \eta_1 $. 

Taking two more $z$ derivatives for the integral \re{K3f2derbefore}, we see that $\pa_{z}^{2} K_{3}f$ is a sum of three terms: \\
\[
   \int_{H^{+}} F (s) \left[   \underbrace{  \int_{1}^{s_{1}} \int_{1}^{\eta_1} \cdots \int_{1}^{\eta_{k-2-\tau_0}} }_\text{$k-1 - \tau_{0}$ integrals} \frac{\pa_{z}^{2} g^{\iota} (z, (\eta_{k-1-\tau_0}, \hat{s} )  ) }{ \wti{\Phi}^{n-l}(z, (\eta_{k-1-\tau_0}, \hat{s} ) )} \, d \eta \right] \, dV(s);
\]
\[
  \int_{H^{+}} F (s) \left[  \underbrace{  \int_{1}^{s_{1}} \int_{1}^{\eta_1} \cdots \int_{1}^{\eta_{k-2-\tau_0}} }_\text{$k-1 - \tau_{0}$ integrals} \frac{ \psi_{1}  \pa_{z} g^{\iota} (z, (\eta_{k-1-\tau_0}, \hat{s} )  )}{\wti{\Phi}^{n-l + 1}(z, (\eta_{k-1-\tau_0}, \hat{s} ) )}  \, d \eta \right] \, dV(s);
\]
\eq{J2}
  \int_{H^{+}} F (s) \left[  \underbrace{  \int_{1}^{s_{1}} \int_{1}^{\eta_1} \cdots \int_{1}^{\eta_{k-2-\tau_0}} }_\text{$k-1 - \tau_{0}$ integrals}  \frac{ \psi_{2} \, g^{\iota} (z, (\eta_{k-1-\tau_0}, \hat{s} )  )}{ \wti{\Phi}^{n-l +2}(z, (\eta_{k-1-\tau_0}, \hat{s} ) ) } \, d \eta \right] \, dV(s), 
\eeq
where $\psi_1$ is a multiple of $\pa_{z} \wti{\Phi} (z, s)$ and $\psi_2 (z, s)$ is a linear combination of $\left( \pa_{z} \wti{\Phi} (z, s) \right)^{2}$ and $\pa_{z}^{2} \wti{\Phi} (z,s) \wti{\Phi} (z, s)$. 
Since $\wti{\Phi} (z,s) = \Phi(z, \phi_{z}^{-1} (s))$, and $\Phi$, $\phi$ (\re{cdphidef}) are holomorphic in $z \in V$, we see that $\psi_{1}$ and $\psi_2$ are smooth functions in $z \in V$. By \re{divder},  \re{gzs}, \re{West}, \re{detests1'}, we obtain 
\eq{balest1}
  | \pa_{z}^{q} g^{\iota} (z, s)| \leq C s_{1}^{-\nu_{0}} |\phi^{-1} (s) - z|^{1- 2l -q - \mu_{0} - \mu_{2} - \nu_{1}}, \quad q = 0, 1, 2. 
\eeq
Since $\wti{\Phi} (z, s) = \Phi(z, \phi^{-1} (s))$, by \re{Phie2} and \re{Phi}, we have
\eq{balest2}
  C | \phi^{-1} (s) - z| \geq |\wti{\Phi} (z, s)| \geq c \left| \phi^{-1} (s) -z \right|^{2}. 
\eeq 
Replacing $s$ in \re{balest1} and \re{balest2} by $(\eta_{k-1-\tau_{0}}, \hat{s})$, we see that in order to estimate $\pa_{z}^{2} K_{3 }f$ it suffices to estimate the integral in \re{J2} for $l=n-1$. i.e. 
\begin{align} \label{K2fint}
  I (z) &= \int_{H^{+}} F (s) J_{k-1-\tau_0} (z, s) \,  dV(s) 
\end{align}
where $J_{k-1 - \tau_0} (z, s) $ is
\eq{Jk-1-tau0}
  \underbrace{  \int_{1}^{s_{1}} \int_{1}^{\eta_1} \cdots \int_{1}^{\eta_{k-2-\tau_0}}}_\text{$k-1 - \tau_{0}$ integrals} \frac{(\eta_{k-1-\tau_0})^{-\nu_{0}}  \wti{\Phi}^{-3}(z, (\eta_{k-1-\tau_0} , \hat{s} ) )    \, d \eta  }{|\phi^{-1} ( (\eta_{k-1-\tau_0}, \hat{s} ) ) - z|^{(2n-3) + \mu_{0} +  \mu_{2} + \nu_{1}}}. 
\eeq
Observe that
\eq{s1seq}
  \eta_{k-1-\tau_0} \geq \eta_{k-2-\tau_0} \geq  \cdots \geq s_{1}. 
\eeq
Thus
\begin{align} \label{Iest1}
  &\left| \phi^{-1} ((\eta_{k-1-\tau_0}, \hat{s})) -z \right|
  \\ \nonumber & \quad \geq c \left\{ | \eta_{k-1-\tau_0}- \phi^{1} (z) | + | \hat{s} - \hat{\phi} (z)| \right\}  \\
  \nonumber & \quad = c \left\{ \eta_{k-1-\tau_0} + \left| \phi^{1} (z) \right| +  | \hat{s} - \hat{\phi} (z)| \right\}, \: \eta_{k-1-\tau_0} >0 , \: \phi^{1} (z) = \rho(z) < 0 
  \\ \nonumber & \quad \geq c \left\{ s_{1} + \left| \phi^{1} (z) \right| +  | \hat{s} - \hat{\phi}(z)| \right\}, \quad \text{by \re{s1seq}}.
  \\ \nonumber &  \quad = c \left\{ | s_{1}  - \phi^{1} (z) | +  | \hat{s}- \hat{\phi} (z)|  \right\} 
  \\ \nonumber  &  \quad \geq c | s - \phi(z)|, \quad s = (s_{1}, \hat{s}) = (s_1, s_2, t_3, \dots, t_{2n}). 
\end{align}
By \re{Phie1} we have 
\begin{align}  \label{Iest2}
  &\left| \wti{\Phi}(z, (\eta_{k-1-\tau_0}, \hat{s} ) ) \right| \\ \nonumber
  & \quad =  \left| \Phi \left( z, \phi^{-1} ( (\eta_{k-1-\tau_0}, \hat{s} ) ) \right) \right|
  \\ \nonumber & \quad \geq c \left( \right. d(z) + \left| \phi^{1} (\phi^{-1} (\eta_{k-1-\tau_0}, \hat{s} ) ) \right| + \left| \phi^{2}(\phi^{-1} (\eta_{k-1-\tau_0}, \hat{s} )) \right| 
  \\ \nonumber & \qquad + \left| \phi' (\phi^{-1} (\eta_{k-1-\tau_0}, \hat{s} ) ) \right|^{2} \left. \right)
   \\ \nonumber & \quad = c \left( \right. d(z)  + \eta_{k-1-\tau_0}+ | s_{2} | + |t|^{2} \left. \right)
   \\ \nonumber & \quad \geq c \left( d(z) + s_{1} + | s_{2}  | + |t|^{2} \right) , \quad \text{by \re{s1seq}}, 
\end{align} 
where $c$ is independent of $z \in V \cap D$. Since
\eq{Iest3}
  | \phi^{-1} ( (\eta_{k-1-\tau_0}, \hat{s} ) ) - z| \geq c \eta_{k-1-\tau_0}, 
\eeq
we have
\[ 
  |\phi^{-1} (s) - z|^{- \mu_{0} -  \mu_{2} - \nu_{1}} \leq C \left( \eta_{k-1-\tau_0} \right)^{- \mu_{0} - \mu_{2} - \nu_{1}}. 
\] 
Using \re{Iest1} and \re{Iest2} we can estimate the integral $J_{k-1-\tau_{0}}(z,s)$ (\re{Jk-1-tau0}) by pulling out $\left| \phi^{-1} ((\eta_{k-1-\tau_0}, \hat{s})) -z \right|^{-(2n-3)}$ and $\wti{\Phi}^{-3}(z, (\eta_{k-1-\tau_0} , \hat{s} ) )$ from the integral sign. In view of \re{Iest3} and  $\nu_{0} + \nu_{1} + \mu_{0} + \mu_{2} = k-1 - \tau_{0}$,  we obtain for $F(s) \in L^{p}(H^{+})$, 
\begin{align*}
   | I (z) | &\leq C(D) \int_{H^{+}} \frac{ | F (s) | \left| s- \phi(z) \right|^{- (2n-3)}  }{(d(z) + s_{1} + |s_{2}| + |t|^{2})^{3}} \, \times 
   \\ & \quad  \left[  \underbrace{  \int_{s_1}^{1} \int_{\eta_1}^{1} \cdots \int_{\eta_{k-2-\tau_0}}^{1} }_\text{$k-1 - \tau_{0}$ integrals} \frac{d \eta}{(\eta_{k-1-\tau_0})^{k-1 - \tau_{0}}} \right] \, dV(s) 
   \\ &\leq  C(D) \int_{H^{+}} |F (s)| A(z, s) ( 1 + | \log s_{1} |) \, dV(s), 
\end{align*}
where we denote
\[ 
   A(z, s) =  \frac{\left| s- \phi(z) \right|^{- (2n-3)} }{(d(z) + s_{1}+ |s_{2}| + |t|^{2})^{3}}, \quad s= (s_{1}, s_{2}, t). 
\]
By the second inequality in \re{Phie2}, we have
\begin{align} \label{s-phiz}
  |s - \phi(z) |  
  &\geq c | \phi^{-1} (s) - z |  \\
  & \nonumber \geq c \{ d(z) + s_1 +  | \phi^{2} (\phi^{-1} (s)) | + | \phi'(\phi^{-1} (s) ) | \} \\
  & \nonumber \geq c \{ d(z) + s_1 + |s_2| +  |t| \}. 
\end{align} 
We estimate 
\[ 
  |I(z)| \leq C(D) \int_{H^{+}}  | F(s) | \left[ A(z, s) \right]^{\frac{1}{p}} \left[ A(z, s) \right]^{\frac{1}{p'}} (1 + \left| \log s_{1} \right|) \, dV(s). 
\] 
By Hölder's inequality, 
\[  
 |I(z)|^{p} \leq  \left[ \int_{H^{+}} A (z, s) | F (s) |^{p} \, ds \right]  \left[ \int_{H^{+}} A(z, s) (1 + | \log s_{1}|)^{p'} ds \right]^{\frac{p}{p'}}. 
\]
Using polar coordinates for $t = (t_{1}, \dots, t_{2n-2})$, and \re{s-phiz}, we have 
\begin{align*}
 &\int_{H^{+}} A(z, s) ( 1 + | \log s_{1}|)^{p'} \, ds
 \\ & \quad \leq \int_{0}^{1} \int_{0}^{1}  \int_{0}^{1} \frac{ (1 + | \log s_{1} |)^{p'} t^{2n-3} \,  ds_{1} \, ds_{2} \, dt}{(d(z) + s_{1} + s_{2} + t^{2})^{3}(d(z) + s_{1} + s_{2} + t)^{2n - 3}}  \\
 & \quad \leq C(\ve)  \int_{0}^{1} \int_{0}^{1}  \int_{0}^{1} \frac{ s_{1}^{-\ve p'} t^{2n-3} \,  ds_{1} \, ds_{2} \, dt }{(d(z) + s_{1} + s_{2} + t^{2})^{3}(d(z) + s_{1} + s_{2} + t)^{2n - 3}}  \\
 & \quad \leq C'(\ve) d(z)^{-\yh - \ve p'}, 
\end{align*}
where for the last inequality we used \rl{intestl} (ii) with $\all = 1 - \ve p'$. 
Thus 
\begin{align} \label{Hqhdsetup} 
&\int_{D} d(z)^{\gm p} |\pa_{z}^{2} Kf(z)|^{p} \, dV(z)  
\\ \nonumber & \qquad \leq C(D) [C(\ve)]^{\frac{p}{p'}} \int_{D} \int_{H^{+}} d(z)^{\gm'}  A(z, s) |F(s)|^{p} \, dV(s) \, dV(z)
\\ \nonumber & \qquad = C(D) [C(\ve)]^{\frac{p}{p'}} \int_{H^{+}} \left[ \int_{D}  d(z)^{\gm'} A(z, s) \, dV(z)   \right] |F(s)|^{p} \, dV(s), 
\end{align}
where
\[
  \gm' = \gm p - \left( \yh \right)  \left( \frac{p}{p'}\right) - \ve p
  = \left( \gm - \yh - \ve \right)p + \yh. 
\]
Pick $\gm$ and $\ve$ with $\gm > \yh + \ve$, then $\gm' > \yh$. For each $s \in H^{+}$, let $\wti{\phi}_{\phi^{-1}(s)} : V \to V$ be the coordinate map given by \re{phitizeta}:
\[ 
 \wti{\phi}_{\phi^{-1}(s)}^{1} (z) = \rho(z), \quad \wti{\phi}_{\phi^{-1}(s)}^{2} (z) = Im (\rho_{\phi^{-1}(s)} \cdot (\phi^{-1} (s) -z)),    
\]
\[ 
   \wti{\phi}_{\phi^{-1}(s)}^{'} (z) = (Re(\phi^{-1}(s) - z), \, Im (\phi^{-1}(s)- z) ).  
\]
Write $\wti{\phi}_{\phi^{-1}(s)} (z) = (\ti{s}_{1}, \ti{s}_{2}, \ti{t})$. Using polar coordinates for $(\ti{s}_{1}, \ti{s}_{2}) \in \R^{2}$, $\ti{t} \in \R^{2n-2}$, and $ c d (z) \leq \ti{s}_{1} \leq C d(z)$, we get for $\gm' > \yh$, 
\begin{align*}
&\int_{D} d(z)^{\gm'} A (z, s) \, dV(z) 
\\ &\qquad \leq C (D) \int_{\ti{t} \in [-1,1]^{2n-2}} \int_{\ti{s}_{2} = -1}^{1} \int_{\ti{s}_{1} = 0}^{1} \frac{ \ti{s_{1}}^{\gm'} d\ti{s}_{1} \, d\ti{s}_{2} \, d  \ti{t}}{(d(\zeta) + \ti{s}_{1} + |\ti{s}_{2}| + \ti{t}^{2})^{n+1}}
\\ &\qquad \leq C(D) \int_{0}^{1} \int_{0}^{1} \frac{\ti{s}^{1+ \gm'} \ti{t}^{2n-3}}{(d(\zeta) + \ti{s} + \ti{t}^{2})^{n+1}} \, d\wti{s} \, d\wti{t} 
\leq C(D, \gm'), 
\end{align*} 
where for the last inequality we used \rl{intestl} (i). We also have
\begin{align*}
  \left[ \int_{D} d(z)^{\gm p} |\pa_{z}^{2} Kf(z)|^{p} \, dV(z) \right]^{\frac{1}{p}} 
  &\leq C(D, \gm) \left[  \int_{U \sm D} | F(s) |^{p}  \, dV(s) \right]^{\frac{1}{p}} \\
  &\leq C(D, \gm) \| \var \|_{W^{k,p} (D)} , 
\end{align*}
Thus we have shown that
\[ 
\| u_{1} \|_{W^{1,p}_{\beta}(D)} \leq C(D, \beta  ) \| \var \| _{W^{1,p}(D)}, \quad \quad \text{for any $0 < \beta < \yh$. }
\qquad  \qed
\]

\section{Estimates for $H_{0}$}
\le{H0ele}
(i) Let $0 < \del < 1$, and $n \geq 2$. Then
\[ 
   \int_{0}^{1} \int_{0}^{1} \frac{s \, t^{2n-3}}{(\del + s + t^{2})^{n+1}} \, dt \, ds \leq C(n) (1 + | \log \del |). 
\]
(ii) Let $\all > 0$, $0 < \del < 1$, and $n \geq 2$. Then 
\[
    \int_{0}^{1} \int_{0}^{1} \frac{s^{1+\all} t^{2n-3}}{(\del + s + t^{2})^{n+1}} \, dt \, ds \leq C(n, \all). 
\]
(iii) Let $0 < \all < 1$.  Then
\[ 
   \int_{0}^{1} \int_{0}^{1} \int_{0}^{1} \frac{s_{1}^{-1+ \all} t^{2n-3} }{(\del + s_{1} + s_{2} + t^{2})^{n+1}} \, ds_{1} \, ds_{2} \, dt
   \leq C(n, \all) \del^{-1 + \all}. 
\] 
\ele
\noindent $\mathit{Proof.}$
(i) Denote the integral by $I$ and split the domain of integration $[0,1] \times [0,1]$ into six regions: \\
$R_{1}: \del \leq t^{2} \leq s$. We have
\[ 
   I \leq \int_{\del}^{1} \int_{t=0}^{\sqrt{s}} \frac{s t^{2n-3}}{s^{n+1}} \, dt \, ds
   \leq C(n) \int_{\del}^{1} s^{-1} \, ds \leq C(n) (1 + | \log \del |). 
\] 
$R_{2}: t^{2} \leq \del \leq s$. We have
\[
I \leq \int_{0}^{\sqrt{\del}} \int_{s= \del}^{1} \frac{s t^{2n-3}}{s^{n+1}} \, ds \,dt
\leq  \int_{0}^{\sqrt{\del}}  \del^{-n+1} t^{2n-3} \, dt \leq C. 
\]
$R_{3}: \del \leq s \leq t^{2} $. We have
\[ 
   I \leq \int_{\sqrt{\del}} ^{1} \int_{s=0}^{t^{2}} \frac{s t^{2n-3}}{t^{2n+2}} \, ds \, dt
   \leq C \int_{\sqrt{\del}} ^{1} t^{-1} \, dt \leq C ( 1+ | \log \del |). 
\]

$R_{4}: s \leq \del \leq t^{2}$. We have
\[ 
  I \leq \int_{0}^{\del} \int_{t=\sqrt{\del}}^{1} \frac{s t^{2n-3}}{t^{2n+2}} \, dt \, ds
  \leq C \int_{0}^{\del} \del^{-2} s \, ds \leq C. 
\]
$R_{5}: t^{2} \leq s \leq \del$. We have
\[ 
  I \leq \int_{0}^{\del} \int_{t=0}^{\sqrt{s}} \frac{s t^{2n-3}}{\del^{n+1}} \, dt \, ds
  \leq C(n) \int_{0}^{\del} \del^{-(n+1)} s^{n} \, ds \leq C (n). 
\]
$R_{6}: s \leq t^{2} \leq \del$. We have
\[ 
   I \leq \int_{0}^{\sqrt{\del}} \int_{s=0}^{t^{2}} \frac{s t^{2n-3}}{\del^{n+1}} \, ds \, dt
   \leq C \int_{0}^{\sqrt{\del}} \del^{-(n+1)} t^{2n +1} \, dt \leq C(n). 
\]
(ii) Split the domain of integration $[0,1] \times [0,1]$ into six regions. \\
$R_{1}: \del \leq t^{2} \leq s$. We have
\[ 
   I \leq \int_{\del}^{1} \int_{t=0}^{\sqrt{s}} \frac{s^{1+\all} t^{2n-3}}{s^{n+1}} \, dt \, ds
   \leq C(n) \int_{\del}^{1} s^{\all-1} \, ds \leq C (n, \all). 
\] 
$R_{2}: \del \leq s \leq t^{2} $. We have
\[ 
I \leq \int_{\sqrt{\del}} ^{1} \int_{s=0}^{t^{2}} \frac{s^{1+ \all} t^{2n-3}}{t^{2n+2}} \, ds \, dt
\leq C \int_{\sqrt{\del}} ^{1} t^{2 \all -1} \, dt \leq C(\all).  
\] 
$R_{3}: t^{2} \leq \del \leq s$. We have
\[
  I \leq \int_{0}^{\sqrt{\del}} \int_{s= \del}^{1} \frac{s^{1+\all} t^{2n-3}}{s^{n+1}} \, ds \,dt
  \leq  \int_{0}^{\sqrt{\del}}  \del^{-n+1 + \all} t^{2n-3} \, dt \leq C(n). 
\]
$R_{4}: s \leq \del \leq t^{2}$. We have 
\[ 
I \leq \int_{0}^{\del} \int_{t=\sqrt{\del}}^{1} \frac{s^{1+ \all} t^{2n-3}}{t^{2n+2}} \, dt \, ds
\leq C \int_{0}^{\del} \del^{-2} s^{1+ \all} \, ds \leq C \del^{\all}. 
\]
$R_{5}: t^{2} \leq s \leq \del$. We have
\[ 
I \leq \int_{0}^{\del} \int_{t=0}^{\sqrt{s}} \frac{s^{1+ \all} t^{2n-3}}{\del^{n+1}} \, dt \, ds
\leq C(n) \int_{0}^{\del} \del^{-(n+1)} s^{n + \all} \, ds \leq C (n) \del^{\all}. 
\]
$R_{6}: s \leq t^{2} \leq \del$. We have
\[ 
I \leq \int_{0}^{\sqrt{\del}} \int_{s=0}^{t^{2}} \frac{s^{1+ \all} t^{2n-3}}{\del^{n+1}} \, ds \, dt
\leq C \int_{0}^{\sqrt{\del}} \del^{-(n+1)} t^{2n + 2 \all + 1} \, dt \leq C(n) \del^{\all}. 
\]
(iii) 
Divide the domain of integration $[0,1] \times [0,1] $ into four regions: \\
$R_{1}: t^{2} > \del, s_{1}, s_{2}$. We have 
\begin{align*}
   I &\leq \int_{\sqrt{\del}}^{1} \frac{t^{2n-3}}{t^{2n+2}} \left( \int_{0}^{t^{2}} s_{1}^{-1 + \all} \, ds_{1} \right) \left( \int_{0}^{t^{2}} \, ds_{2} \right) \, dt  \\
   &\leq \int_{\sqrt{\del}}^{1}  t^{-3 + 2 \all} \, dt 
   \leq C(\all) \del^{-1 + \all}. 
\end{align*}
$R_{2}: \del > t^{2}, s_{1}, s_{2}$. We have
\begin{align*}
  I &\leq \del^{-(n+1)} \int_{0}^{\sqrt{\del}} t^{2n-3} \left( \int_{0}^{\del} s_{1}^{-1 + \all} \, ds_{1} \right) \left( \int_{0}^{\del} \, ds_{2} \right) \, dt 
  \\ &\leq C(n, \all) \del^{-(n+1) + n-1 + \all + 1 } 
  \leq C(n, \all) \del^{-1 + \all}. 
\end{align*}
$R_{3}: s_{1} > \del, t^{2}, s_{2}$. We have
\begin{align*}
   I &\leq \int_{\del}^{1} \frac{s_{1}^{-1 + \all}}{s_{1}^{n+1}} \left( \int_{0}^{\sqrt{s_{1}}} t^{2n-3} \, dt \right) \left( \int_{0}^{s_{1}} \, ds_{2} \right)  \, ds_{1} \\
   &\leq C \int_{\del}^{1} s_{1}^{\all -2} 
   \leq C \del^{-1+ \all}.    
\end{align*}
$R_{4}: s_{2} > \del, t^{2}, s_{1}$. We have 
\begin{align*}
 I &\leq \int_{\del}^{1} s_{2}^{-(n+1)} \left( \int_{0}^{\sqrt{s_{2}}} t^{2n-3} \, dt \right) \left( \int_{0}^{s_{2}} s_{1}^{-1 + \all} \, ds_{1} \right) \, ds_{2} 
 \\ &\leq C(\all) \int_{\del}^{1} s_{2}^{\all -2} \, ds_{2} 
 \leq C(\all) \del^{-1 + \all}. \qquad \qquad \qquad  \qquad \qquad \qquad \qed
\end{align*}
We now prove the estimate for the holomorphic projection operator $H_{0}$. In this case we have a loss which is arbitrarily small in the expoenent of the weight. 
\th{H0thm} 
  Let $ D \subset \subset \C^{n}$ be a bounded strictly pseudoconvex domain with $C^2$ boundary. Let $H_{0} \var$ be defined by formula \re{H0}. \\
  (i) For any $1 < p < \infty$, we have
  \[ 
     \| H_{0} \var \|_{W^{0,p}_{\beta}(D)} \leq C(D, \beta) \| \var \|_{W^{1,p} (D)}, \quad \text{for any $\beta$, \: $0 < \beta <1$.} 
  \]  
  (ii) Suppose $2n< p < \infty$, and $k \geq 2$ is an integer. We have
  \[ 
      \| H_{0} \var \|_{W^{k-1,p}_{\beta}(D)} \leq C (D, \beta) \| \var \|_{W^{k,p} (D)}, \quad \text{for any $\beta$, \: $0 < \beta <1$.} 
  \] 
\eth
\begin{proof}
(i) 
In view of \re{Om1}, $H_{0} \var$ can be written as a linear combination of 
\[ 
  Kf (z) = \int_{U \sm D} f(\zeta) \frac{\hat{\pa}_{\zeta} W (z,\zeta)}{\Phi^{n} (z, \zeta) } \, dV(\zeta), \quad \Phi(z, \zeta) = W \cdot (\zeta -z), 
\]
where $f$ denotes a coefficient function of $[\dbar, E] \var$. Thus $f \in L^{p}(U \sm D)$, and $f \equiv 0$ in $D$. Let  $\hat{\pa}_{\zeta} W$ denote the products of $W$ and its first derivatives in  $\zeta$. Let $W_{1} = (\hat{\pa}_{\zeta} W, \pa_{z}^{\nu_{0}} \hat{\pa}_{\zeta} W)$. 

Let $V$ be a neighborhood of $\zeta_{\ast} \in bD$, as given by the remark after \rp{regLeray}. Using a partition of unity in $\zeta$ and $z$ space, we can assume
\[
  \text{supp}_{\zeta} \, f \subset V \sm D, \quad \quad \text{supp}_{z} \, \Om^{01}_{0,q} (z, \zeta) \subset V \cap D. 
\] 
We have
\begin{align} \label{subintk}
&\int_{D} d(z)^{\gm p} | \pa_{z} Kf (z) \var (z)|^{p} \, dV(z)
\\ \nonumber  &\quad =  \int_{D} d(z)^{\gm p} \left| \int_{U \sm D} f(\zeta) k(z, \zeta) d V(\zeta)  \right|^{p} \, d V(z) 
\\ \nonumber&\quad \leq  \int_{D} d(z)^{\gm p} \left| \int_{U \sm D} |f(\zeta)| |k(z, \zeta)|^{\frac{1}{p}} |k(z, \zeta)|^{\frac{1}{p'}} d V(\zeta)  \right|^{p} \, d V(z)
\\ \nonumber &\leq \int_{D} d(z)^{\gm p} \left[ \int_{U \sm D} | f(\zeta) |^{p} |k(z, \zeta)| \, dV(\zeta) \right] \left[  \int_{U \sm D} | k(z, \zeta)| \, dV(\zeta) \right]^{\frac{p}{p'}} \, dV(z), 
\end{align}
where we set
\[ 
  k(z, \zeta) = \pa_{z} \left(  \frac{\hat{\pa}_{\zeta} W (z,\zeta)}{\Phi^{n} (z, \zeta) } \right)
  = \frac{\pa_{z} \hat{\pa}_{\zeta} W(z, \zeta)}{\Phi^{n} (z, \zeta)} - n \frac{\hat{\pa}_{\zeta} W(z, \zeta) \pa_{z}  \Phi(z, \zeta) }{\Phi^{n+1}(z, \zeta)}. 
\]
For fixed $z \in V \cap D$, define the coordinate map $\phi_{z}: V \to \phi(V)$ as in \re{cdphidef}. Write $\phi_{z} (\zeta) = (s_{1}, s_{2},t)$. Then from \re{West} and \re{Philbe}, we have 
\eq{H0kest}
  | k(z, \zeta) | \leq \frac{C(D)}{\left( d(z)+ s_{1} + |s_{2} | + |t|^{2} \right)^{n+1}}. 
\eeq
Integrating using polar coordinates for $s = (s_{1}, s_{2}) \in \R^{2}$ and $t=(t_{1}, \dots, t_{2n-2}) \in \R^{2n-2}$, we have by \rl{H0ele} (i),
\begin{align*}  
  \int_{U \sm D} | k(z, \zeta)| \, dV(\zeta) 
  &\leq C(D) \int_{t \in [-1,1]^{2n-2}} \int_{-1}^{1} \int_{0}^{1} \frac{ ds_{1} ds_{2} \, dt}{(d(z) + s_{1} + |s_{2}| + t^{2})^{n+1}}  \\
  &\leq C(D) \int_{0}^{1} \int_{0}^{1} \frac{s \, t^{2n-3}}{(d(z) + s + t^{2})^{n+1}} \, ds \, dt \\  
  &\leq C(D) (1 + | \log d(z)|) \\
  &\leq C(D, \ve) d(z)^{-\ve}, 
\end{align*} 
for any $\ve>0$. Substituting the above estimate into the last line of \re{subintk} we get 
\begin{align*}  
&\int_{D} d(z)^{\gm p} | \pa_{z} H_{0} \var (z)|^{p} \, dV(z)
\\ & \qquad \leq [C(D, \ve)]^{\frac{p}{p'}} \int_{D} d(z)^{\gm'}  \left[ \int_{U \sm D} | f(\zeta) |^{p} |k(z, \zeta)| \, dV(\zeta) \right] \, dV(z) 
\\ & \qquad = [C(D, \ve)]^{\frac{p}{p'}} \int_{U \sm D}  | f(\zeta) |^{p} \left[ \int_{D} d(z)^{\gm'} |k(z, \zeta)| \, dV (z) \right] \, dV(\zeta), 
\end{align*} 
where we set 
\begin{gather} \label{H0fdgm'}
   \gm' = \gm p  - \ve \frac{p}{p'} = p (\gm - \ve) + \ve. 
\end{gather}
Choose $\gm$ and $\ve$ with $\gm > \ve > 0$ so that $\gm' > 0$. For $\zeta \in V \sm D$, let $\wti{\phi}_{\zeta}$ be the coordinate map given by \re{phitizeta}. $\wti{\phi}_{\zeta} (z) = (\ti{s}_1, \ti{s}_2, \ti{t})$. Using polar coordinates for $\ti{s}= ( \ti{s}_{1}, \ti{s}_{2}) \in \R^{2}$, $ \ti{t} \in \R^{2n-2}$, and $ c d (z) \leq \ti{s}_{1} \leq C d(z)$, we get for $\gm' > 0$, 
\begin{align} \label{lastest}
  &\int_{D} d(z)^{\gm'} |k (z, \zeta) | \, dV(z)
  \\ \nonumber & \qquad \leq C (D) \int_{\ti{t} \in [-1,1]^{2n-2}} \int_{\ti{s}_{2} = -1}^{1} \int_{\ti{s}_{1} = 0}^{1} \frac{ \ti{s}_{1}^{\gm'} d\ti{s}_{1} \, d\ti{s}_{2} \, dt}{(d(z) + \ti{s}_{1} + |\ti{s}_{2}| + \ti{t}^{2})^{n+1}} 
  \\ \nonumber & \qquad \leq C(D) \int_{0}^{1} \int_{0}^{1} \frac{ \ti{s}^{1+ \gm'} \ti{t}^{2n-3}}{(d(\zeta) + \ti{s} + | \ti{t}| ^{2})^{n+1}} \, d \ti{s} \, d \ti{t} 
  \leq C(D, \gm')
\end{align}
where in the last inequality we used \rl{H0ele} (ii). Consequently 
\[ 
  \left[ \int_{D} d(z)^{\gm p} | \pa_{z} H_{0} \var (z)|^{p} \, dV(z) \right]^{\frac{1}{p}} \leq C(D, \gm) 
  \left[ \int_{U \sm D} | f(\zeta) |^{p} \, dV(\zeta) \right]^{\frac{1}{p}}, 
\]
i.e. 
\[  \| H_{0} \var \|_{W^{0,p}_{\beta} (D)} \leq C(D, \beta) \| \var \|_{W^{1,p} (D)}, \] 
for any $\beta$, $0 < \beta < 1$.
\\ \\
(ii) Assume 
\[
\text{supp}_{\zeta} \, f \subset V \sm D, \quad \quad \text{supp}_{z} \, \Om^{01}_{0,q} (z, \zeta) \subset V \cap D, 
\]  
where $V$ is the same as in (i). We can write $\pa_{z}^{k} H_{0} \var$ as a linear combination of 
\[ 
  K_{1} f (z) =  \int_{U \sm \ov{D}} f (\zeta) \frac{W_{1} (z, \zeta)}{\Phi^{n+ l}(z, \zeta)}  \, dV(\zeta), \quad 0 \leq l \leq k, 
\] 
where $W_{1}(z, \zeta)$ denotes some polynomial in $\pa_{z}^{k_{0}} \hat{\pa}_{\zeta} W(z, \zeta)$ and $\pa_{z}^{k_{1}} \Phi(z, \zeta)$, for $k_0, k_1 \geq 0$. 

We now integrate by parts to reduce the exponent of $\Phi$ in the denominator to $n+1$. 
Let $\zeta_{i_{\ast}}$ be such that $u(z, \zeta) := \pa_{\zeta_{i_{\ast}}} \Phi (z, \zeta) \neq 0$ for $z \in V \cap D$ and $\zeta \in V \sm D$. Write
\[ 
   K_{1} f(z) = \int_{U \sm \ov{D}} f(\zeta) W_{1} (z, \zeta) u^{-1}(z, \zeta) \pa_{\zeta_{i_{\ast}}} \Phi^{-(n+l -1) } (z, \zeta) \, dV(\zeta). 
\] 
By \rp{Sobemb}, $f \in W^{k-1,p} (U) \subset C^{k-2+ \all} (U)$, for $\all = 1- \frac{2n}{p} \in(0,1)$. Since $f \equiv 0$ in $\ov{D}$, we have $|f(\zeta)| \leq | f |_{U; k-2+\all} d(\zeta)^{k-2+ \all}$, for $\zeta \in U \sm \ov{D}$. Here $d(\zeta) = \dist(\zeta, D)$. By \re{West}, $|\pa_{\zeta}^{i} W_{1}(z, \zeta)| \leq C(D) d(\zeta)^{-i}$, and $|\pa_{\zeta}^{i} u^{-1} (z, \zeta)| \leq C(D) d(\zeta)^{-i}$, for $i=0,1,2,\dots$. In particular, 
\[ 
  \left| W_{1} u^{-1} \right| \leq C(D), \quad
  \left| \pa_{\zeta_{i_{\ast}}} \left(  W_{1} u^{-1}  \right) \right| \leq C(D) d(\zeta)^{-1}.   
\]
In view of \re{Phie1} for fixed $z \in D$, we have $| \Phi^{-(n+k-1)} (z, \zeta)| \leq C (z)$ and 
\[ 
  \left| \pa_{\zeta_{i_{\ast}}} \Phi^{-(n+l-1)} \right| = \left| \Phi^{-(n+l)} u \right| \leq  C(z, D). 
\]
Thus $W_{1}u^{-1}$ and $\Phi^{-(n+ l-1)}$ satisfy the estimate \re{gbloword} for $i_{k} = 0$. Applying \rl{ibple} (i) and (ii) we obtain
\begin{align*}
  K_{1} f(z) &= \int_{U \sm D} f(\zeta) \pa_{\zeta_{i_{\ast}}} \left( W_{1} (z, \zeta) u^{-1} (z, \zeta) \Phi^{-(n+l -1) } (z, \zeta) \right) dV(\zeta) 
  \\ &\quad - \int_{U \sm D} f(\zeta) \pa_{\zeta_{i_{\ast}}} \left( W_{1}  (z, \zeta) u^{-1}  (z, \zeta) \right) \Phi^{-(n+l-1)} (z, \zeta) \, dV(\zeta) 
  \\
  &= - \int_{U \sm D} \left( \pa_{\zeta_{i_{\ast}}} f (\zeta)  \right) W_{1} (z, \zeta) u^{-1} (z, \zeta)  \Phi^{-(n+l -1) } (z, \zeta)  \, dV(\zeta) \\
  &\quad - \int_{U \sm D}  f (\zeta)   \pa_{\zeta_{i_{\ast}}} \left( W_{1}  (z, \zeta) u^{-1}  (z, \zeta) \right)  \Phi^{-(n+l -1) } (z, \zeta)  \, dV(\zeta). 
\end{align*}

We can repeat this procedure $l - 1 (\leq k-1)$ times. Indeed, suppose we have done $m$ times, $1 \leq m \leq k-2$. Then the integral is a linear combination of
\begin{align*}
&\int_{U \sm D} \left( \pa_{\zeta_{i_{\ast}}}^{m_{1}} f \right) \pa_{\zeta_{i_{\ast}}}^{m_{2}} \left\{ u^{-1}, W_{1} \right\} \Phi^{-(n+l - m )} \, dV(\zeta), \quad \quad m_{1} + m_{2} = m, 
\\ \nonumber &\qquad =  \int_{U \sm D} \left( \pa_{\zeta_{i_{\ast}}}^{m_{1}} f \right) \pa_{\zeta_{i_{\ast}}}^{m_{2}} \left\{ u^{-1}, W_{1} \right\} u^{-1} \pa_{\zeta_{i_{\ast}}}  \Phi^{-(n+l - m - 1)} \, dV(\zeta), 
\end{align*}
where $ \pa_{\zeta_{i_{\ast}}}^{m_{2}} \left\{ u^{-1}, W_{1} \right\}$ is a linear combination of 
\begin{gather} \nonumber
\pa_{\zeta_{i_{\ast}}}^{\la_1} (u^{-1}) \pa_{\zeta_{i_{\ast}}}^{\la_1} (u^{-1}) \cdots \pa_{\zeta_{i_{\ast}}}^{\la_{p}} (u^{-1})
\pa_{\zeta_{i_{\ast}}}^{\la_0} (W_{1} (z, \zeta)), 
\\ \nonumber
\la_{i} \geq 0, \quad \quad \sum_{i=0}^{p} \la_i = m_2 . 
\end{gather} 
We have $\pa_{\zeta_{i_{\ast}}}^{m_{1}} f \in W^{1,p}(U \sm D) \cap C^{k-2 - m_{1} +\all} (\ov{U \sm D})$, and $\pa_{\zeta_{i_{\ast}}}^{m_{2}} \left\{ u^{-1}, W_{1} \right\}$ satisfies the estimate \re{gbloword} for $i_k = m_{2}$, and $u^{-1}$, 
$\Phi^{-(n+l - m - 1)}$ satisfy estimates \re{gbloword} for $i_k = 0$. Since $k-2- m_{1} - m_{2} = k-2- m \geq 0$, the hypothesis of \rl{ibple} (i) and (ii) hold, and we can do the procedure one more time. 
 In the end we can write $K_{1} f$ as a linear combination of 
\[
   \wti{K_{2}} f(z) = \int_{U \sm \ov{D}} [\pa_{\zeta_{i_{\ast}}}^{k-1} f (\zeta)] W_{1} ( u^{-1})^{k-1} \Phi^{-(n+1)} (z, \zeta) \, dV(\zeta),  
\]    
and 
\begin{gather} \label{H0K2f} 
  K_{2} f(z) = \int_{U \sm \ov{D}} \pa_{\zeta_{i_{\ast}}}^{\tau_{0}} f(\zeta) \pa_{\zeta_{i_{\ast}}}^{\tau_{1}} \left\{ u^{-1}, W_{1} \right\}  \Phi^{-(n+1)} (z, \zeta) \, dV(\zeta), 
\\ \nonumber
  \tau_{0} + \tau_{1} = l -1 \leq k-1. 
\end{gather}
As $ \pa_{\zeta_{i_{\ast}}}^{k-1} f \in L^{p}(U)$, $\wti{K_{2}} f$ can be estimated in the same way as part (i):
\begin{align} \label{H0K2tie}
   \left[ \int_{D} d(z)^{\gm p} | \wti{K_{2}} f (z)|^{p} \, dV(z) \right]^{\frac{1}{p}} 
   & \leq C (\gm, D) \| f \|_{W^{k-1, p} (U \sm D)} \\
   \nonumber &\leq C (\gm, D) \| \var \|_{W^{k,p} (D)}
\end{align} 
for any $\gm > 0$. For $K_{2} f$, we integrate by parts in the direction $s_{1}$. Take $V$ and $\phi$ as in \re{cdphidef}, and set $\hat{\phi} = (\phi^{2}, \dots, \phi^{2n}) $. Let $U_0 = V \cap (U \sm D)$. Define
\[ 
  \wti{\pa_{\zeta_{i_{\ast}}}^{\tau_{0}} f }(s) = \pa_{\zeta_{i_{\ast}}}^{\tau_{0}} f (\phi^{-1} (s)), \quad
  \wti{\Phi} (z, s) = \Phi(z, \phi^{-1} (s)), 
\] 
\[ 
  g(z,s) = \pa_{\zeta_{i_{\ast}}}^{\tau_{1}} \left\{ u^{-1}, W_{1} \right\} (z, \phi^{-1} (s)) \left| det(D \phi^{-1}) (s) \right|. 
\]
Then we have, for $s= (s_{1}, \hat{s}), \; \hat{s} = (s_{2}, t_3, \dots, t_{2n})$,
\begin{align} \label{K2h0exp}
   K_{2} f(z) &= \int_{H^{+}}  \wti{\pa_{\zeta_{i_{\ast}}}^{\tau_{0}} f }(s) g(z,s) \wti{\Phi}^{-(n+1)}(z, s) \, dV(s) \\
   \nonumber &= \int_{H^{+}} \wti{\pa_{\zeta_{i_{\ast}}}^{\tau_{0}} f }(s) \pa_{s_{1}} I_{1} (z, s) \, dV(s), 
\end{align}
where 
\[
  I_{1}(z, s) = \int_{1}^{s_{1}} g (z, (\eta_1, \hat{s})  ) \wti{\Phi}^{-(n+1)} (z, (\eta_1, \hat{s}) ) \, d\eta_1. 
\]
By \re{West} and \re{detests1'}, we have
\eq{H0gest}
  \left| g(z, (\eta_1, \hat{s})) \right| \leq C(D) \eta_{1}^{- \tau_{1}}. 
\eeq
Thus
\begin{align*}
   | I_{1}(z,s) |
    \leq C(D, z) \int_{s_{1}}^{1} \eta_{1}^{-\tau_{1}} \, d \eta_1
    \leq 
    \begin{cases}
    C(D, z) s_{1}^{-(\tau_{1} -1)} &  \text{if $\tau_{1} > 1$}, \\ 
    C(D, z) | \log s_{1} | &  \text{if $\tau_{1} = 1$}. \\
    C(D, z) & \text{if $\tau_1 = 0$,}
    \end{cases}
\end{align*}  
and 
\[ 
  | \pa_{s_{1}} I_{1} (z, s) | = \frac{|g (z, (\eta_1, \hat{s}) )| }{\left| \wti{\Phi}^{n+1} (z, (\eta_1, \hat{s}) ) \right| }  
  \leq 
   \begin{cases}
  C(D, z) s_{1}^{- \tau_{1} } &  \text{if $\tau_{1} > 1$}, \\ 
  C(D, z) s_{1}^{-1} & \text{if $\tau_{1} = 1$}. \\
  C(D, z) & \text{if $\tau_1 = 0$.}
  \end{cases}
\] 
If $\tau_{1} >1$, we have $\wti{\pa_{\zeta_{i_{\ast}}}^{\tau_{0}} f } \in W^{1,p} (H^{+}) \cap C^{j+\all}(H^{+}) $ for $j = k-2 - \tau_{0} \geq 0$ ($\tau_{0} \leq k-2$) and $I_{1}$ satisfies estimate \re{gbloword} with $i_k$ replaced by $\tau_{1} -1$. Furthermore, we have
\[ 
  j - (\tau_{1} -1) = k-1 - \tau_{0} - \tau_{1} \geq k -1 - l \geq 0.
\]
Thus we can apply \rl{ibple} (ii) and integrate by parts in \re{K2h0exp} to get
\eq{K2h0ibpafter1}
  K_{2} f (z) = - \int_{H^{+}} \pa_{s_{1}} \wti{\pa_{\zeta_{i_{\ast}}}^{\tau_{0}} f}(s) 
  I_{1} (z, s) \, dV (s). 
\eeq
If $\tau_{1} = 1$, then $I_{1}$ satisfies estimate \re{ibplog}. Thus we can apply \rl{ibple} (iii) and integrate by parts to get \re{K2h0ibpafter1}. If $\tau_0 = 0$, then $I_1$ satisfies estimate \re{gbloword} with $i_k$ being replaced by $0$ and again we can integrate by parts by \rl{ibple} (ii) to obtain \re{K2h0ibpafter1}. 

We can integrate by parts $k-1- \tau_{0} = \tau_{1}$ times. Indeed, suppose we have done it $m$ times, $1 \leq m \leq k-2- \tau_0$. We have
\begin{align*}
  K_2 f &= \pm \int_{H^{+}} \pa_{s_{1}}^{m} \wti{\pa_{\zeta_{i_{\ast}}}^{\tau_{0}} f}(s) I_{m} (z,s) \, dV(s) \\
  &=  \pm \int_{H^{+}} \pa_{s_{1}}^{m} \wti{\pa_{\zeta_{i_{\ast}}}^{\tau_{0}} f}(s) \pa_{s_{1}} I_{m+1} (z,s) \, dV(s),
\end{align*}
where
\[
I_{m} (z, s) = \underbrace{  \int_{1}^{s_{1}} \int_{1}^{\eta_1} \cdots \int_{1}^{\eta_{m-1}} }_\text{$m$ integrals}  \frac{g (z, ( \eta_m,  \hat{s}  ) ) \, [d \eta]^{m} }{\wti{\Phi}^{n+1}(z, (\eta_{m}, \hat{s} ) ) }, 
\]
and we denote $[d \eta]^{m} := \, d \eta_m \cdots d \eta_1$. By \re{H0gest} and \re{Phie1},  
\begin{align*}
  |I_{m+1} (z,s) | &\leq  C (D,z) \underbrace{  \int_{s_{1}}^{1} \int_{\eta_1}^{1} \cdots \int_{\eta_{m}}^{1} }_\text{$m+1$ integrals}   
  \eta_{m+1}^{- \tau_{1}} \, [d \eta]^{m+1} \\
  &\leq     
  \begin{cases}
  C(D, z) s_{1 }^{- \tau_{1} + (m+1)}, & \text{if $m+1 < \tau_{1}$, }   \\
  C(D,z) ( 1+ | \log s_{1} |), & \text{if $m+1 = \tau_{1}$},  \\
  C(D,z), & \text{if $m+1 > \tau_{1}$}, 
  \end{cases}
\end{align*}
and
\begin{align*}
 | \pa_{s_{1}} I_{m+1} (z,s) | = | I_{m} (z,s) |
 \leq \begin{cases}
 C(D, z) s_{1 }^{- \tau_{1} + m}, & \text{if $m+1 < \tau_{1} $},   \\
 C(D,z) s_{1}^{-1}, & \text{if $m+1 = \tau_{1}$}, \\
 C(D,z) ( 1+| \log s_{1}|), & \text{if $m+1 > \tau_{1}$}. 
 \end{cases}
\end{align*}
Applying \rl{ibple} (ii) and (iii) to these cases we obtain 
\[ 
K_{2} f (z) = \pm \int_{H^{+}}  \left( \pa_{s_{1}}^{m+1} \wti{\pa^{\tau_{0}}_{\zeta_{i}} f}  \right) (s) I_{m+1} (z, s) dV(s). 
\] 
In conclusion, we can integrate by part $k-1-\tau_0$ times to transform $K_2 f $ (\re{H0K2f}) to the form 
\gan
  K_{2} f (z) = - \int_{H^{+}} F (s) I_{k-1-\tau_0} (z, s) \, dV(s),  
\end{gather*}
where 
\gan
  F (s) = \pa_{s_{1}}^{k-1 - \tau_{0}} \wti{ \pa_{\zeta_{i_{\ast}}}^{\tau_{0}}} f (z, s)  \in L^{p} (U \sm D),
\end{gather*}
and 
\gan
  I_{k-1-\tau_0} =  \underbrace{  \int_{1}^{s_{1}} \int_{1}^{\eta_{1}} \cdots \int_{1}^{\eta_{k-2-\tau_0}} }_\text{$k-1 - \tau_{0}$ integrals} \frac{g (z, (\eta_{k-1-\tau_{0}}, \hat{s}) )}{ \wti{\Phi}^{n+1}(z, (\eta_{k-1-\tau_0}, \hat{s})) } \, \left[ d \eta \right]  
\end{gather*}
with $[ d \eta] = d\eta_{k-1- \tau_0} \, \cdots \, d\eta_1$. By \re{Iest2}, we have
\eq{H0phihatest}
  \left| \wti{\Phi} (z, (\eta_{k-1-\tau_0}, \hat{s})) \right| \geq c \left( d(z) + s_{1} + |s_{2}| + |t|^{2} \right),  
\eeq
where $c$ is independent of $z \in V$. From \re{H0gest} and \re{H0phihatest} we see that $  | K_{2} f(z) |$ is bounded by
\gan
   C(D) \int_{H^{+}}  \frac{| F(s) | \, dV(s)}{ \left( d(z) + s_{1} + |s_{2} | + |t|^{2} \right)^{n+1} } \left[  \underbrace{  \int_{s_1}^{1} \int_{\eta_1}^{1} \cdots \int_{\eta_{k-2-\tau_0}}^{1} }_\text{$k-1-\tau_0$ integrals}  (\eta_{k-1-\tau_0})^{- \tau_{1}} \, d \eta  \right] 
   \\ \leq  C(D)  \int_{H^{+}}  | F(s)| \left[ k(z, s) \right] \,  ( 1 + | \log s_{1} | )  \, dV( s ), \quad \tau_1 \leq k -1 - \tau_0, 
\end{gather*}
where we denote
\[ 
   k(z, s) = \frac{1}{ \left( d(z) + s_{1} + |s_{2}| + |t| ^{2} \right)^{n+1} }. 
\]
Using polar coordinates for $t = (t_{1}, \dots, t_{2n-2})$, and applying \rl{H0ele} (iii) for $\all = 1 - \ve p'$, we have for any $\ve > 0$, 
\begin{align*}
   \int_{H^{+}} [k(z, s)] (1 + | \log s_{1} |) ^{p'} \, dV(s) 
   &\leq C(\ve)  \int_{0}^{1} \int_{0}^{1} \int_{0}^{1} \frac{s_{1}^{- \ve p'} \, t^{2n-3}  \, ds_{1}  \, ds_{2} \, dt  }{(d(z) + s_{1} + s_{2} + t^{2})^{n+1}} \\
   &\leq C(\ve) d(z)^{-\ve p'}. 
\end{align*}
Thus for any $\gm > 0$, we have
\begin{align*}
  &\int_{D} d(z)^{\gm p} | K_{2} f(z) |^{p} \, dV(z) 
  \\ &\quad \leq C(D) \int_{D} d(z)^{\gm p} \left[  \int_{H^{+}} | F(s) | [k(z, s)] ( 1 + | \log s_{1}|) \, dV(s) \right]^{p} \, dV(z)  \\
  &\quad = C(D) \int_{D} d(z)^{\gm p} \left[ \int_{H^{+}} |F(s)| [k(z, s)]^{\frac{1}{p}} [k(z, s)]^{\frac{1}{p'}} 
  (1 + | \log s_{1} |) dV(s) \right]^{p} \, dV(z) \\
  &\quad \leq C(D) \int_{D} d(z)^{\gm p} \left[ \int_{H^{+}} | F(s) |^{p} [k(z, s)] \, dV(s) \right] 
  \\ & \qquad \qquad \times \left[  \int_{H^{+}} k(z, s) ( 1 + | \log s_{1} |)^{p'} \, dV(s) \right]^{\frac{p}{p'}} \, dV(z) \\
  &\quad \leq C(D) [C(\ve)]^{\frac{p}{p'}} \int_{D}  d(z)^{\gm'} \left[ \int_{H^{+}} | F(s) |^{p} [k(z, s)] \, dV(s) \right] \, dV(z) \\
  &\quad \leq C(D) [C(\ve)]^{\frac{p}{p'}}  \int_{H^{+}} |F(s) |^{p} \left[ \int_{D} d(z)^{\gm'} k(z , s) \, dV(z) \right] \, dV(s), 
\end{align*}
where we denote 
\eq{H0hdestgm'}
  \gm' = \gm p - \ve p = p (\gm - \ve). 
\eeq 
Choose $\gm$ and $\ve$ such that $\gm > \ve$. Then $\gm' > 0$.  For each $s \in H^{+}$, let $\wti{\phi}_{\phi^{-1}(s)} : V \to V$ be the coordinate map given by \re{phitizeta}, and $\wti{\phi}_{\phi^{-1}(s)} (z) = (\ti{s}_{1}, \ti{s}_{2}, \ti{t})$. Using polar coordinates for $(\ti{s}_{1}, \ti{s}_{2}) \in \R^{2}$, $\ti{t} \in \R^{2n-2}$, and $ c d (z) \leq s_{1} (z) \leq C d(z)$, we get
\begin{align*} 
   \int_{D} d(z)^{\gm'} k(z , s) \, dV(z)  
   &\leq C \int_{\ti{t} \in [-1, 1]^{2n-2} }  \int_{\ti{s}_{2}=-1}^{1} \int_{\ti{s}_{1} = 0}^{1} \frac{\ti{s}_{1}^{\gm'} \, d\ti{s}_{1} \, d\ti{s}_2 \, d \ti{t} } {\left( d(\zeta) + \ti{s}_{1} + |\ti{s}_{2}| + |\ti{t}|^{2} \right)^{n+1}}  \\
   &\leq C \int_{0}^{1} \int_{0}^{1} \frac{\ti{s}^{1+  \gm'} \ti{t}^{2n-3}} {\left( d(\zeta) + \ti{s} + \ti{t}^{2} \right)^{n+1}} \, d \ti{s} \, d \ti{t}  \\
   &\leq C(n, \gm'), 
\end{align*} 
where in the last inequality we used \rl{H0ele} (ii). Hence
\[ 
   \left[ \int_{D} d(z)^{\gm p} | K_{2} f(z) |^{p} \, dV(z) \right]^{\frac{1}{p}}  \leq C(D, \gm) \left[ \int_{U \sm D} |F(\zeta) |^{p} \, dV(\zeta) \right]^{\frac{1}{p}}. 
\]
Combine this and \re{H0K2tie} we obtain for any $\gm > 0$, 
\[ 
   \left[ \int_{D} d(z)^{\gm p} | \pa_{z}^{k} H_{0} \var(z) |^{p} \, dV(z) \right]^{\frac{1}{p}} \leq C(D, \gm) \left[ \int_{U \sm D} |F(\zeta) |^{p} \, dV(\zeta) \right]^{\frac{1}{p}}, 
\] 
or
\[
   \| H_{0} \var \|_{W^{k-1,p}_{\beta} (D)} \leq C(D, \beta) \| \var \|_{W^{k,p} (D)}, \quad 0 < \beta < 1.
\]
\end{proof} 
\newcommand{\doi}[1]{\href{http://dx.doi.org/#1}{#1}}
\newcommand{\arxiv}[1]{\href{https://arxiv.org/pdf/#1}{arXiv:#1}}

\def\MR#1{\relax\ifhmode\unskip\spacefactor3000 \space\fi%
\href{http://www.ams.org/mathscinet-getitem?mr=#1}{MR#1}}

\bibliographystyle{plain}

\begin{bibdiv}
	\begin{biblist}
		
		\bib{A-S79}{article}{
			author={Ahern, Patrick},
			author={Schneider, Robert},
			title={Holomorphic {L}ipschitz functions in pseudoconvex domains},
			date={1979},
			ISSN={0002-9327},
			journal={Amer. J. Math.},
			volume={101},
			number={3},
			pages={543\ndash 565},
			url={https://doi.org/10.2307/2373797},
			review={\MR{533190}},
		}
		
		\bib{AW06}{article}{
			author={Alexandre, William},
			title={{$C^k$}-estimates for the {$\overline\partial$}-equation on2
				convex domains of finite type},
			date={2006},
			ISSN={0025-5874},
			journal={Math. Z.},
			volume={252},
			number={3},
			pages={473\ndash 496},
			url={https://doi.org/10.1007/s00209-005-0812-y},
			review={\MR{2207755}},
		}
		
		\bib{CD89}{article}{
			author={Chang, Der-Chen~E.},
			title={Optimal {$L^p$} and {H}\"{o}lder estimates for the {K}ohn
				solution of the {$\overline\partial$}-equation on strongly pseudoconvex
				domains},
			date={1989},
			ISSN={0002-9947},
			journal={Trans. Amer. Math. Soc.},
			volume={315},
			number={1},
			pages={273\ndash 304},
			url={https://doi.org/10.2307/2001384},
			review={\MR{937241}},
		}
		
		\bib{C-S01}{book}{
			author={Chen, So-Chin},
			author={Shaw, Mei-Chi},
			title={Partial differential equations in several complex variables},
			series={AMS/IP Studies in Advanced Mathematics},
			publisher={American Mathematical Society, Providence, RI; International
				Press, Boston, MA},
			date={2001},
			volume={19},
			ISBN={0-8218-1062-6},
			review={\MR{1800297}},
		}
		
		\bib{G-T01}{book}{
			author={Gilbarg, David},
			author={Trudinger, Neil~S.},
			title={Elliptic partial differential equations of second order},
			series={Classics in Mathematics},
			publisher={Springer-Verlag, Berlin},
			date={2001},
			ISBN={3-540-41160-7},
			note={Reprint of the 1998 edition},
			review={\MR{1814364}},
		}
		
		\bib{GX18}{article}{
			author={Gong, Xianghong},
			title={H\"{o}lder estimates for homotopy operators on strictly
				pseudoconvex domains with {$C^2$} boundary},
			date={2019},
			ISSN={0025-5831},
			journal={Math. Ann.},
			volume={374},
			number={1-2},
			pages={841\ndash 880},
			url={https://doi.org/10.1007/s00208-018-1693-9},
			review={\MR{3961327}},
		}
		
		\bib{G-S77}{book}{
			author={Greiner, P.~C.},
			author={Stein, Elias~M.},
			title={Estimates for the {$\overline \partial $}-{N}eumann problem},
			publisher={Princeton University Press, Princeton, N.J.},
			date={1977},
			ISBN={0-691-08013-5},
			note={Mathematical Notes, No. 19},
			review={\MR{0499319}},
		}
		
		\bib{HT91}{article}{
			author={Horiuchi, Toshio},
			title={The imbedding theorems for weighted {S}obolev spaces. {II}},
			date={1991},
			ISSN={0579-3068},
			journal={Bull. Fac. Sci. Ibaraki Univ. Ser. A},
			number={23},
			pages={11\ndash 37},
			url={https://doi-org.ezproxy.library.wisc.edu/10.5036/bfsiu1968.23.11},
			review={\MR{1118940}},
		}
		
		\bib{LG09}{book}{
			author={Leoni, Giovanni},
			title={A first course in {S}obolev spaces},
			series={Graduate Studies in Mathematics},
			publisher={American Mathematical Society, Providence, RI},
			date={2009},
			volume={105},
			ISBN={978-0-8218-4768-8},
			url={https://doi.org/10.1090/gsm/105},
			review={\MR{2527916}},
		}
		
		\bib{L-R80}{article}{
			author={Lieb, Ingo},
			author={Range, Michael~R.},
			title={L\"{o}sungsoperatoren f\"{u}r den Cauchy-Riemann-Komplex
				mit $C^{k}$ {A}bsch\"{a}tzungen},
			date={1980},
			ISSN={0025-5831},
			journal={Math. Ann.},
			volume={253},
			number={2},
			pages={145 \ndash 164},
			url={https://doi-org.ezproxy.library.wisc.edu/10.1007/BF01578911},
			review={\MR{597825}}, 
		}
	
		\bib{MJ91}{article}{
			author={Michel, Joachim},
			title={Integral representations on weakly pseudoconvex domains},
			date={1991},
			ISSN={0025-5874},
			journal={Math. Z.},
			volume={208},
			number={3},
			pages={437\ndash 462},
			url={https://doi.org/10.1007/BF02571538},
			review={\MR{1134587}},
		}
		
		\bib{M-P90}{article}{
			author={Michel, Joachim},
			author={Perotti, Alessandro},
			title={{$C^k$}-regularity for the {$\overline\partial$}-equation on
				strictly pseudoconvex domains with piecewise smooth boundaries},
			date={1990},
			ISSN={0025-5874},
			journal={Math. Z.},
			volume={203},
			number={3},
			pages={415\ndash 427},
			url={https://doi.org/10.1007/BF02570747},
			review={\MR{1038709}},
		}
		
		\bib{M-S99}{article}{
			author={Michel, Joachim},
			author={Shaw, Mei-Chi},
			title={A decomposition problem on weakly pseudoconvex domains},
			date={1999},
			ISSN={0025-5874},
			journal={Math. Z.},
			volume={230},
			number={1},
			pages={1\ndash 19},
			url={https://doi.org/10.1007/PL00004685},
			review={\MR{1671846}},
		}
		
		\bib{PK91}{article}{
			author={Peters, Klaus},
			title={Solution operators for the {$\overline\partial$}-equation on
				nontransversal intersections of strictly pseudoconvex domains},
			date={1991},
			ISSN={0025-5831},
			journal={Math. Ann.},
			volume={291},
			number={4},
			pages={617\ndash 641},
			url={https://doi.org/10.1007/BF01445231},
			review={\MR{1135535}},
		}
		
		\bib{RM92}{article}{
			author={Range, Michael~R.},
			title={On {H}\"{o}lder and {BMO} estimates for {$\overline \partial$} on
				convex domains in {$\bf C^2$}},
			date={1992},
			ISSN={1050-6926},
			journal={J. Geom. Anal.},
			volume={2},
			number={6},
			pages={575\ndash 584},
			url={https://doi.org/10.1007/BF02921578},
			review={\MR{1189045}},
		}
		
		\bib{H-R71}{article}{
			author={Romanov, A.~V.},
			author={Henkin, G.~M.},
			title={Exact {H}\"{o}lder estimates of the solutions of the {$\bar
					\delta $}-equation},
			date={1971},
			ISSN={0373-2436},
			journal={Izv. Akad. Nauk SSSR Ser. Mat.},
			volume={35},
			pages={1171\ndash 1183},
			review={\MR{0293121}},
		}
		
		\bib{SY74}{article}{
			author={Siu, Yum-Tong},
			title={The {$\bar \partial $} problem with uniform bounds on
				derivatives},
			date={1974},
			ISSN={0025-5831},
			journal={Math. Ann.},
			volume={207},
			pages={163\ndash 176},
			url={https://doi-org.ezproxy.library.wisc.edu/10.1007/BF01362154},
			review={\MR{0330515}},
		}
		
		\bib{SE70}{book}{
			author={Stein, Elias~M.},
			title={Singular integrals and differentiability properties of
				functions},
			series={Princeton Mathematical Series, No. 30},
			publisher={Princeton University Press, Princeton, N.J.},
			date={1970},
			review={\MR{0290095}},
		}
		
		\bib{WS89}{article}{
			author={Webster, Sidney~M.},
			title={A new proof of the {N}ewlander-{N}irenberg theorem},
			date={1989},
			ISSN={0025-5874},
			journal={Math. Z.},
			volume={201},
			number={3},
			pages={303\ndash 316},
			url={https://doi.org/10.1007/BF01214897},
			review={\MR{999729}},
		}
		
	\end{biblist}
\end{bibdiv}


\end{document}